\documentclass[english,12pt]{amsart}

\usepackage{amsmath,amsfonts,bbm,amsthm}
\usepackage[latin1]{inputenc}  
\usepackage[T1]{fontenc}       
\usepackage[psamsfonts]{amssymb}
\usepackage{dsfont}

\usepackage{latexsym}
\usepackage[dvips]{epsfig}
\usepackage{dcolumn}
\usepackage{tabularx}
\usepackage{longtable}
\usepackage{graphics}
\usepackage{pstricks}
\usepackage{pst-node}
\usepackage{latexsym}
\usepackage{babel}
\usepackage{euscript}

\usepackage{tikz} 
\usetikzlibrary{arrows.meta,calc}

\renewcommand{\theequation}{\thesection.\arabic{equation}}
\newtheorem{theorem}{Theorem}
\newtheorem{lemma}{Lemma}
\newtheorem{proposition}{Proposition}

\newtheorem{remark}{Remark}

\newtheorem{theoA}{Theorem}

\newtheorem{propA}[theoA]{Proposition}

\newtheorem{lemA}[theoA]{Lemma}

\newcommand{\eqnsection}{
\renewcommand{\theequation}{\thesection.\arabic{equation}}
    \makeatletter
    \csname  @addtoreset\endcsname{equation}{section}
    \makeatother}
\eqnsection

\def\r{{\mathbb R}}

\newcommand{\ind}{\mathds{1}}

\def\sss{{\mathcal S}}
\def\demi{{1\over 2}}

\def\Z{{{\Bbb Z}}}

\def\R{{{\Bbb R}}}
\def\C{{{\Bbb C}}}
\def\S{{\mathbb S}}

\def\E{{{\Bbb E}}}

\def\indic{{{\mathbbm 1}}}
\def\Id{{\hbox{Id}}}

\def\cE{{\mathcal E}}

\newtheorem{conjecture}{Conjecture}

\newtheorem{corol}{Corollary}[section]
\newtheorem{rema}{Remark}[section]
\newtheorem{exa}{Example}

\renewcommand{\le}{\leqslant}
\renewcommand{\leq}{\leqslant}
\renewcommand{\ge}{\geqslant}
\renewcommand{\geq}{\geqslant}
\renewcommand{\subset}{\subseteq}

\newcommand{\tW}{\tilde{W}}

\newcommand{\bal}{\begin{align*}}
\newcommand{\eal}{\end{align*}}
\newcommand{\beq}{\begin{eqnarray*}}
\newcommand{\eeq}{\end{eqnarray*}}
\newcommand{\bte}{\begin{theorem}}
\newcommand{\ete}{\end{theorem}}
\newcommand{\bl}{\begin{lemma}}
\newcommand{\el}{\end{lemma}}
\newcommand{\bd}{\begin{description}}
\newcommand{\ed}{\end{description}}

\newcommand{\bc}{\begin{cases}}
\newcommand{\ec}{\end{cases}}
\newcommand{\bp}{\begin{proof}}
\newcommand{\ep}{\end{proof}}
\newcommand{\bco}{\begin{corol}}
\newcommand{\eco}{\end{corol}}

\newcommand{\T}{\mathbb{T}}

\def\bdes{\begin{description}}
\def\edes{\end{description}}

\def\ggg{{\mathcal{G}}}
\def\l{\left}
\def\r{\right}

\newcommand{\gb}{\beta}

\def\beqn{\begin{eqnarray}}
\def\eeqn{\end{eqnarray}}
\def\H{{\mathbb{H}}}
\def\cV{{\mathcal{V}}}
\def\cF{{\mathcal{F}}}
\def\i{{\imath}}
\begin{document}
\author[C. Sabot]{Christophe Sabot}
\address{{Institut Camille Jordan, CNRS UMR 5208}, {Universit\'e Claude Bernard Lyon 1}, {{43 Boulevard du 11 novembre 1918}, {Villeurbanne}, {69622}, {France}}} 
\address{Institut Universitaire de France}
\email{sabot@math.univ-lyon1.fr}
\title[$\H^{2|2}$ and Anderson with Cauchy potential]
{ The $\H^{2|2}$-model and the Anderson model with Cauchy potential: similarities and a conjectured relation.}
\maketitle
\begin{abstract}

We highlight strange similarities between the Anderson model with Cauchy potential with energy-disorder parameters $(E,\lambda)$ and the $\H^{2|2}$ supersymmetric sigma field with initial non local conductances proportional to $W_{x,y}=\lambda^2\l\vert \l(-\Delta-E+\i\lambda)\r)_{x,y}^{-1}\r\vert^2$  where 
$\Delta$ is the discrete Laplace operator. More precisely, looking at the Brownian interpolation of both, we show that they follow very similar equations. Motivated by that we investigate the different phases of the $\H^{2|2}$-model on $\Z^d$ with conductances $rW$ for some $r>0$ and show that whatever $r>0$ the  $\H^{2|2}$-model displays localized and delocalized phases which are coherent with what is expected or proved for the Cauchy potential with energy-disorder parameters $(E,\lambda)$. 
Finally, we propose a conjecture that compares key quantities of the two models and would imply existence of delocalized regime for the Anderson model with Cauchy potential in dimension $d\ge 3$.
\end{abstract}
\section*{Introduction}
The aim of this paper is to highlight some similarities between the Anderson model with Cauchy potential (also called Lloyd's model)  and the  $\H^{2|2}$ supersymmetric sigma field. The Anderson model is a fundamental model of mathematical physics which describes the displacement of a quantum particle in a disordered media. It has attracted a great deal of attention among mathematicians and physicists. The localized phase is now fairly well understood for smooth distribution of potential, the first multidimensional results in this direction were based on the famous multi-scale method of Fr\"olich and Spencer (\cite{FS}), and later the fractional moment method of Aizenman and Molchanov (\cite{Aizenman_Molchanov,Aizenman94}). On the other hand the  existence of a delocalized phase was a major open problem solved only very recently by Openai for the uniform potential (\cite{openai}). See e.g. the monographies on the topic \cite{Aizenman_Warzel} or \cite{Kirsch}.

On the other hand, the $\H^{2|2}$ supersymmetric sigma field was introduced by Zirnbauer, and first analyzed by Disertori Spencer and Zirnbauer (\cite{zirnbauer91,DSZ10,DS10}), as a toy model for the family of supersymmetric sigma models which are expected to describe random band matrices (see in particular the discussion in \cite{SZ4} and \cite{Efetov,DSZ10}).  The $\H^{2|2}$ sigma model is however not formally connected to the Anderson model or random band matrices, while it is expected to be phenomenologically similar. In particular, a remarkable result of Disertori Spencer and Zirnbauer exibits a delocalized phase in dimension $d\ge 3$ at weak disorder, and Disertori Spencer (\cite{DS10}) proved exponential localization in any dimension at strong disorder. In \cite{ST15}, Sabot and Tarr\`es proved that the $\H^{2|2}$-model is related to a reinforced process, the Vertex Reinforced Jump Process (VRJP) and also to the Edge Reinforced Random Walk introduced by Coppersmith and Diaconis. This connection motivated several developments. Particularly useful here will be the representation of the $\H^{2|2}$-model in terms of a random Schr\"odinger operator (with one-dependent potential) at the ground state (\cite{STZ17,SZ19}).

\subsection*{Brownian interpolation of the Anderson model with Cauchy potential} Let us explain informally the similarities between the Anderson model with Cauchy potential (we just say Anderson model in the sequel but it always restricts to the Cauchy potential) and the $\H^{2|2}$-model and the main results we prove. Consider $\ggg=(\cV,\cE)$ a finite non-directed graph and $\Delta$ the discrete Laplacian on $\ggg$ (we are typically interested in a large box of $\Z^d$),  and $V=(V_x)_{x\in \cV}$ the random potential with i.i.d. centered Cauchy random variables with scale ${1\over \sqrt{2}}$. For $\lambda>0$ a disorder parameter and $E\in \R$ the spectral parameter, we consider the Schr\"odinger operator
\begin{eqnarray}\label{GV0}
H_V=V+{1\over \lambda} (-\Delta-E).
\end{eqnarray}
(It is more convenient in our context to consider $H$ as above instead of the more standard  $-\Delta-E+\lambda V$). 
As it is well known the standard Cauchy distribution can be represented as the law of the first hitting point of the real line by a complex Brownian motion starting from $\i\over \sqrt{2}$ where $\i$ is the complex unit. More precisely, for all $x\in \cV$, we consider
$$
X_x(t)=\exp(B_x(t)-t/2),
$$
the geometric Brownian motion, where $(B_x(t))_{x\in \cV}$ are independent standard Brownian motions, and $Y_x(t)$ defined by $Y_x(0)=0$ and 
$$
dY_x(t)=X_x(t) d\tilde B_x(t),
$$
where $(\tilde B_x(t))_{x\in \cV}$ are independent standard Brownian motions, independent of $(B_x)$. Then, we set, with $\i$ the complex unit,
$$
V_x(t)={1\over \sqrt{2}} (Y_x(t)+\i X_x(t)),
$$
which is a complex martingale with square bracket $\left< d \overline{V_x}(t), dV_x\right> =X_x^2(t) dt$. Then $\lim_{t\to \infty} V_x(t):=V_x\in \R$ exists and is a Cauchy random variable with scale ${1\over \sqrt{2}}$. We thus consider a dynamical version of the Anderson model with Cauchy potential and set
$$
W^A_{x,y}(t)=\vert G_{V(t)}(x,y)\vert^2,
$$
where $G_{V(t)}= (H_{V(t)})^{-1}=(V(t)+{1\over \lambda} (-\Delta-E))^{-1}$ is the Green function of $H_{V(t)}$.

By Itô calculus, it is easy to see that $(W(t),X(t))$ satisfies the following equation
\begin{eqnarray}\label{eq_W_A_intro}
\begin{cases}
dW_{x,y} (t)=(W(t) X(t)^2W(t))_{x,y} dt +dM_{x,y}(t)
\\
dX(t)=X^2(t) dB(t),
\end{cases}
\end{eqnarray}
where $X(t)^2$ above is the operator of multiplication by $X(t)^2$ and $dM_{x,y}(t)$ is a local martingale term expressed in terms of $G_{V(t)}$ and $dV(t)$, see Section~\ref{s_SDE_A} (so $M_{x,y} (t)$ is a priori not measurable in terms of $W(t),X(t)$ and $(X(t),W(t))$ is not a diffusion process). We also consider the following expression
$$
\tilde W_{x,y}(t)=X_x(t) W_{x,y}(t) X_y(t),
$$
which satisfies the following equation for all moment $p>0$,
\begin{eqnarray}\label{eq_tildeW_A}
d\tW_{i,j}^p&=& p^2 \tW_{i,j}^{p-1} (\tilde W^2)_{i,j}dt  - p^2  \tilde W_{i,j}^{p} \left((\tW\cdot 1)_i+(\tW\cdot 1)_j\right)dt 
\\
\nonumber 
&&+ p(p-1)\tilde W_{i,j}^p dt+ p\tilde W_{i,j}^p(dB_i(t)+dB_j(t)) +p^2 \delta_{i=j} \tW_{i,j}^p dt \\
\nonumber
&&+ \tW_{i,j}^{p-1} d\tilde M_{i,j}(t)
\end{eqnarray}
where $(\tW\cdot 1)_i=\sum_j W_{i,j}$ and $d\tilde M_{i,j}(t)=X_i(t) X_j(t) dM_{i,j}(t)$ is a local martingale term.

\subsection*{Brownian interpolation of the $\H^{2|2}$ model} In \cite{STZ17,SZ19}, the $\H^{2|2}$ model was shown to be equivalent to a random Schr\"odinger operator. We will not define the $\H^{2|2}$ model but rather focus on its random Schr\"odinger representation. More precisely, consider now a finite set $\cV$ and a symmetric matrix with non-negative coefficients $W=(W_{x,y})_{x,y\in \cV}$. For $\beta=(\beta_x)_{x\in \cV}$ a potential on $\cV$ we define
$$
H_\beta= H^W_\beta= \beta-W,\;\;\; \hbox{ and } \;\;\; G_\beta=G_\beta^W=(H_\beta)^{-1}
$$
its Green function when it is well-defined. We write $H_\beta >0$ when the operator is positive definite, in which case $G_\beta$ has non negative coefficients.
The following measure defined on potentials $\beta=(\beta_x)_{x\in \cV}$ on $\cV$ is a probability distribution (see \cite{STZ17}, Definition~1)
\begin{eqnarray}\label{nu_W}
\nu_\cV^{W}(d\beta)= \l(\frac{2}{\pi} \right)^{|\cV|/2}
{\ind_{H_\beta>0}
\over 
{\sqrt{|H_\beta|}}}
\exp\l(-\frac{1}{2}\langle 1 ,H_\beta 1 \rangle \r) d\beta
\end{eqnarray}
where $1$ above is the vector of $1$'s, so that $\langle 1 ,H_\beta 1 \rangle =\sum_i \beta_i-\sum_i\sum_j W_{i,j}$.

On a large box of the lattice $\Z^d$ and $W_{x,y}=\indic_{x\sim y} W$  for a real $W>0$, the results of Disertori Spencer and Zirnbauer (\cite{DSZ10,DS10}) implies the existence of phase transition for $G_\beta$ in dimension $d\ge 3$, between a localized phase with exponential decay of $G_\beta$ for $W$ large enough and a delocalized phase for $W$ small enough, see later for precise definitions. (Note that different proofs of the localized phase were provided either using reinforced processes, \cite{ACK14}, or the fractional moment method \cite{Collevecchio_Zeng21}).

The similarities with the Anderson model appear most clearly through the Brownian motion representation of the random Schr\"odinger operator $H_\beta$ introduced in \cite{SZ20,GRSZ24}. We need first to fix some notations. For a vector of "times" $T=(T_x)_{x\in \cV}\in \R_+^\cV$, we define 
$$
K_T=K^W_T= (\Id-TW),
$$
where $T$ represents the operator of multiplication by $(T_x)$ (i.e. the diagonal matrix with coefficients $(T_x)$). We then define the process $X(t)=(X_x(t))_{x\in \cV}$ by the following SDE
$$
dX_x(t)=X_x(t) dB_x(t)-X^2_x(t) (W(t)X(t))_x dt,
$$
where
$$
W(t)= W (K_{T(t)}^W)^{-1}, 
$$
and for all $x\in \cV$,
$$
T_x(t)= \int_{0}^t (X_x(s))^2 ds,
$$
is the square bracket of $X_x(t)$. Hence $W(0)=W$. The equation is well defined as shown in Theorem~2.1 of \cite{GRSZ24}, in particular for all time $t>0$, $T(t)$ is such that $K^W_{T(t)}>0$, so there is no explosion. The connection with the probability distribution defined above, \eqref{nu_W}, is that $1/T=(1/T_x)_{x\in \cV}$ has distribution $\nu_\cV^{W}$ where $T_x:=\lim_{t\to\infty} T_x(t)$ is the total square bracket of $X_x(t)$. Behind the proof of this statement, there is the fact that the process $(X(t),W(t))$ 
is a Markov process, following the equation below
\begin{eqnarray}\label{eq_W_S}
\begin{cases}
dW_{x,y} (t)=(W(t) X(t)^2W(t))_{x,y} dt
\\
dX_x(t)=X_x^2(t) dB(t)-X_x^2(t)(W(t)X(t))_x dt.
\end{cases}
\end{eqnarray}
The equation above is rather different from \eqref{eq_W_A_intro}, however working with these equations one can see several computational coincidences. In particular $(W(t)X(t))_x$ is a martingale for both. It can be seen either by It\^o formula or for the Anderson model using the Ward identity $(W(t)X(t))_x=-\sqrt{2}\Im(G_{V(t)}(x,x))$. In both cases, the asymptotic of this martingale contains key information on the localization or delocalization of the models.

But the formal similarities appear more clearly on the quantity $\tilde W_{x,y}(t)$ defined as before by
$$
\tilde W_{x,y}(t)=X_x(t) W_{x,y}(t) X_y(t).
$$
Note that this transformation is natural since it reduces the scaling parameter to $1$ (see below, Section~\ref{sec_Gauge_S}). Then,  $(\tilde W_{x,y}(t))_{x,y\in \cV}$ becomes Markov and its moments satisfy the following SDE
\begin{eqnarray}\label{eq_tildeW_S}
d\tW_{i,j}^p&=& p \tW_{i,j}^{p-1} (\tilde W\cdot \tilde W)_{i,j}dt  - p  \tilde W_{i,j}^{p} \left((\tW1)_i+(\tW1)_j\right)dt 
\\
\nonumber 
&&+ p(p-1)\tilde W_{i,j}^p dt+ p\tilde W_{i,j}^p(dB_i(t)+dB_j(t)) +p^2\delta_{i=j} \tW_{i,i}^p dt.
\end{eqnarray}

\subsection*{Initial value for the $\H^{2|2}$} Compared to \eqref{eq_tildeW_A}, the first line of the right-hand side of \eqref{eq_tildeW_S} has a factor $p$ instead of $p^2$, the second line is the same, and there is no extra martingale term as in the third line of \eqref{eq_tildeW_A}. Note that the change from $p$ to $p^2$ in  \eqref{eq_tildeW_A} comes from the quadratic variation of the martingale terms. Besides, since the first line is $p+1$ homogeneous in $\tilde W$, while the second line is $p$-homogeneous, then the process $\hat W:={1\over p} \tilde W$ satisfies the SDE given exactly by the first two lines of \eqref{eq_tildeW_A}. It naturally encourages us to start the two processes with the same initial condition, which would mean here to take the conductances of the $\H^{2|2}$-model as 
\begin{eqnarray}\label{W_initial}
W^S_{x,y}(0):= r \vert G_{{\i\over \sqrt{2}}}(x,y)\vert ^2= r\l \vert \l({\i\over \sqrt{2}} \Id +{1\over \lambda} (-\Delta-E)\r )^{-1}_{x,y}\r\vert^2
\end{eqnarray}
for some positive real $r$. The equation for the $p$-moment suggests $r={p}$ but there is no obvious reason to take one $p$ rather than another. 
We comment on that in Section~\ref{sec_discussion} and Section~\ref{sec_one_point} and propose two conjectures comparing key quantities of Anderson and the $\H^{2|2}$. Working on these models, we can see several coincidences even beyond the relations between \eqref{eq_tildeW_A} and \eqref{eq_tildeW_S}. However, we could not establish general rigorous connexions.

\subsection*{What we prove: the phase diagram of the corresponding $\H^{2|2}$-model} Using the tools developed to analyse the $\H^{2|2}$-model and the VRJP, it is possible to prove that on $\Z^d$, $d\ge 3$, with non-local conductances given by \eqref{W_initial} (hence the conductances are positive for all pairs of vertices $x$ and $y$ but have exponential decay), the $\H^{2|2}$-model has a delocalized phase and a localized phase coherent with what is expected for the Anderson model. More precisely, we prove that for any choice of $r>0$, for any $E\in \mathring{\sigma(\Delta)}=(0,4d)$ (for technical reasons we remove $d$ points to $(0,4d)$), there exists $\lambda_0(E,r)$ small enough such that for $\lambda<\lambda_0(E,r)$, the $\H^{2|2}$-model with conductances equal to \eqref{W_initial} is delocalized (precise meaning is given in Section~\ref{sec_loc_deloc}). On the other hand, for all $r>0$ there is a domain containing the points $(E,\lambda)$ with $E\not\in \sigma(-\Delta)$ and $\lambda$ small or $E+{\i\over \sqrt{2}} \lambda$ far enough from $\sigma(-\Delta)$,  on which the $\H^{2|2}$-model is localized (see Theorem~\ref{thm_delocalization} 
and figure~\ref{figure1} 
below).

\subsection*{Important quantities of the models} If we look at the limit values of $W^{A}(t)$ and $W^S(t)$ 
we have
$W^A_{x,y} (\infty)=\l | G_V(x,y) \r |^2$ on the Anderson side and, for $x\neq y$, $W_{x,y}^S(\infty)=(W^S(0) (K_{T(\infty)}^{W^S(0)})^{-1})_{x,y}=({1\over \beta} G^{W^S(0)}_\beta {1\over \beta})_{x,y}$ on the $\H^{2|2}$ side since $T(\infty)=1/\beta$. The expectation of fractional moments of both quantities give key information on the localization of the two models.  

Besides, it is interesting to look at the martingales $(W(t)X(t))_x$ on both sides.  In particular, to avoid degeneracies, we can look at a version where $X(t)$ is frozen on a "boundary" subset of vertices. (Typically, we will use a uniform pinning, which on the Anderson side amounts to considering the Schr\"odinger operator $G_{V+\i\epsilon}$.)  In that case, the limit value $(W^A(\infty ) X^A(\infty))_x$ is given by $-\sqrt{2} \Im(G_{V}(x,x))$, where $V$ is Cauchy distributed on active vertices and $V$ is equal to ${\i\over \sqrt{2}}$ on frozen vertices,
 which can give information on the spectral measure at $x$. On the $\H^{2|2}$ side, $(W^S(\infty ) X^S(\infty))_x$ is given by $(W^S(0) \psi)_x$ where $\psi$ is a martingale of the model which characterizes the recurrence or transience of the VRJP (see Section~\ref{sec_martingale} for details). Hence, we propose the following heuristic correspondence between the Cauchy Anderson model and the $\H^{2|2}$ model  with the specific choice of conductances given by \eqref{W_initial}:
\begin{eqnarray*}
\hbox{Anderson}\;\;\; & & \;\;\; \H^{2|2}
\\
\l | G_V(x,y) \r |^2 \;\;\; &\longleftrightarrow& \;\;\;  ({1\over \beta} G_\beta {1\over \beta})_{x,y}
\\
\l | \sqrt{2} \Im(G_V(x,x)) \r|  \;\;\; &\longleftrightarrow& \;\;\;  (W \psi)_x 
\eeq

In section~\ref{sec_discussion}, we conjecture that for $0<p\le \demi $ there is a real $r_0(p)$, $0<r_0(p)<1$, such that for any vertex transitive finite graph with a uniform pinning $\epsilon>0$, and for any $(E,\lambda)$, there is an $r$ in $(r_0(p),1)$ such that $\E\l(| \sqrt{2} \Im(G_{V+\i\epsilon}(x,x))|^p\r)$, for $x$ a vertex, is equal to the corresponding value for the $\H^{2|2}$ model with conductances $W^S(0)=rW^A(0)$. Together with the delocalization result for the $\H^{2|2}$ model, that would imply existence of absolutely continuous spectrum for the Anderson model with Cauchy potential. In Section~5, we conjecture an explicite value for the optimal lower bound $r_0(p)$, corresponding to the case of one-point with pinning.

\subsection*{Organisation of the paper} In Section~\ref{s_SDE_A} and \ref{s_SDE_S}, we derive the Brownian interpolations of respectively the Anderson model with Cauchy potential and the $\H^{2|2}$-model. In Section~\ref{sec_delocalization}  we state and prove the localization and delocalization results of the $\H^{2|2}$ model on $\Z^d$ with conductances \eqref{W_initial}, and also recall useful   properties of the $\beta$-potential. Section~\ref{sec_discussion} and \ref{sec_one_point} concern discussion and conjectures.

\section{The Brownian interpolation of the Anderson model with Cauchy potential}\label{s_SDE_A}
\subsection{Definitions and equation of $W_{x,y}$}
Let us first fix the context. We consider a discrete finite graph $\ggg=(\cV,\cE)$. We denote by
$\Delta$ the associated discrete Laplacian acting on functions on $\cV$ by
$$
\Delta f(x)=\sum_{y\sim x} (f(y)-f(x)).
$$
(Note that we could add some positive weights 
on the edges.)
We fix $\lambda >0$ a disorder parameter and $E\in \R$ an energy parameter. For $$
V=(V_x)_{x\in \cV} \in \C^\cV,
$$ 
a potential, we set
\begin{eqnarray}\label{GV}
G_{V,E,\lambda}(x,y)=G_V(x,y):= \left( V +{1\over \lambda} (-\Delta-E)\right)_{x,y}^{-1},
\end{eqnarray}
the associated Green function when it is defined.
We are interested in the Anderson model with Cauchy potential, more precisely when 
$V$ is an i.i.d. sequence of random variables with centered Cauchy distribution with scale ${1\over \sqrt{2}}$. As is well-known, the Cauchy random variable with scale ${1\over \sqrt{2}}$ is the law of the point where a complex Brownian motion starting at ${\i\over \sqrt{2}}$ first hits the real line. 
From this Brownian interpolation of the Cauchy law we get an interpolation of the Green function and consider the equation of its square modulus.

 We denote $B(t)=(B_x(t))_{x\in \cV}$ an i.i.d. family of standard Brownian Motions (BM) and $ \tilde B(t)=(\tilde B_x(t))_{x\in \cV}$ another i.i.d. family of standard BM independent of $B(t)$.  We set for $x\in \cV$,
$$
X_x(t):= \exp(B_x(t)- \demi t); \;\;\;  A_x(t)= \int_0^t X_x^2(s) ds .
$$
and
$$
V_x(t) = {1\over \sqrt{ 2}}( Y_x(t) + \i X_x(t)), \;\hbox{ where } \; Y_x(t)=\tilde B_x(A_x^{-1}(t)),
$$
so that $V_x(t)$ is a complex martingale with square bracket 
$$
\left< d\overline V_x(t), dV_x(t)\right >= X_x^2(t) dt,
$$
and $V_x:=\lim_{t\to \infty} V_x(t)={1\over \sqrt{2}} \lim_{t\to \infty} Y_x(t),
$
is a centered Cauchy distribution with scale ${1\over \sqrt{2}}$.
With the notation \eqref{GV} above, we set for all $x,y$ in $\cV$,
$$
 W_{x,y} (t):= \vert G_{V(t)}(x,y) \vert^2.
$$
We will repeatedly use the Ward identity:
\begin{eqnarray}\label{Ward}
-\sqrt{2} \Im(G_{V(t)} (x,x))=(W(t)X(t))_{x}.
\end{eqnarray}
$G_{V(t)} $ is a complex local martingale and we have
$$
dG_{V(t)} (x,y)=- \sum_{z\in \cV} G_{V(t)}(x,z) dV_z ( t) G_{V(t)} (z,y),
$$
and
\begin{eqnarray*}
&&dW_{x,y}(t)\\
 &=& -\sum_{z} G_{V(t)}(x,z) dV_z ( t) G_{V(t)} (z,y) \overline{G_{V(t)} (x,y) }  -\sum_{z}  \overline{G_{V(t)}(x,z) dV_z ( t) G_{V(t)} (z,y)} G_{V(t)} (x,y)
\\
&&+ \sum_{z} W_{x,z}(t) W_{z,y}(t) \left< d\overline V_z(t), dV_z(t)\right>.
\end{eqnarray*}
Finally, we get,
\begin{eqnarray}\label{eq_W_A}
\begin{cases}
dW_{x,y}(t)  = (WX^2W)_{x,y} dt+ dM_{x,y},
\\
dX_x(t)=X_x(t) dB_x(t).
\end{cases}
\end{eqnarray}
where on the right-hand side $X^2$ stands for the operator of multiplication by $X^2(t)$ (i.e. the diagonal matrix with coefficients $(X^2_x(t))$, and
with 
$$
 dM_{i,j}(t) =\sum_k \sigma_{i,j}^k dV_k(t) +\overline{\sigma}_{i,j}^k d\overline{V}_k(t).
 $$ 
and
\begin{eqnarray}\label{sijk}
\sigma_{i,j}^k =  -G_{V(t)}(i,k)G_{V(t)} (k,j) \overline{G_{V(t)} (i,j) }.
\end{eqnarray}

\subsection{The Gauge transform $\tilde W_{x,y}$}\label{sec_Gauge_A}
For all $i,j\in \cV$, we set
$$
\tilde W_{i,j}(t)=X_i(t) W_{i,j}(t) X_j(t).
$$
From the definition of $W(t)$, we see that if we include $X_i(t)$ and $X_j(t)$ in the Green function it changes the potential $V(t)$ to $V(t)/ X(t)={\i\over\sqrt{2}}+Y(t)/X(t)$ and changes the Laplacian term from $-\Delta-E$ to ${1\over \sqrt{X(t)}}(-\Delta-E){1\over \sqrt{X(t)}}$. Hence it brings back the potential to a unit imaginary part, as at time $0$. A similar transformation is very natural on the $\H^{2|2}$ model see Proposition~\ref{prop_Gauge}. 

\begin{proposition}\label{prop_eq_tW_A}
The process $\tilde W(t)$ is solution of the following equation: for all $i, j\in \cV$,
\begin{eqnarray}\label{eq_tW_A}
d\tilde W_{i,j} 
 &=&(\tilde W^2)_{i,j}dt  -\tilde W_{i,j} \left((\tW\cdot 1)_i+(\tW\cdot 1)_j\right)dt 
 \\
\nonumber  &&+ \tilde W_{i,j}(dB_i(t)+dB_j(t))+ \delta_{i,j} W_{ii} dt+ X_iX_j dM_{ij}(t)
\end{eqnarray}
where $(\tW\cdot 1)_i=\sum_j \tW_{i,j}$. More generally, for $p>0$, we have
\begin{eqnarray}\label{eq_tW_A_p}
d\tW_{i,j}^p&=& p^2 \tW_{i,j}^{p-1} (\tilde W^2 )_{i,j}dt  - p^2  \tilde W_{i,j}^{p} \left((\tW\cdot 1)_i+(\tW\cdot 1)_j\right)dt 
\\
\nonumber &&+ p(p-1)\tilde W_{i,j}^p+ p\tilde W_{i,j}^p(dB_i(t)+dB_j(t)) +p^2\delta_{i,j} W_{ii}^p dt
 \\
\nonumber 
&&+ p\tW_{i,j}^{p-1} X_iX_j dM_{i,j}
\end{eqnarray}
where we have as before
\begin{eqnarray}\label{dM}
dM_{i,j}= \sum_{k} \left( \sigma^{ij}_k dV_k(t)+\overline{\sigma^{ij}_k dV_k(t)}\right).
\end{eqnarray}
\end{proposition}
\begin{remark}\label{rk_diag}
Note that in \eqref{eq_tW_A_p} the term $\tilde W_{i,i} \tilde W_{i,j}^p$ in $\tW_{i,j}^{p-1} (\tilde W^2)_{i,j}$ cancels the  term $\tilde W_{i,i} \tilde W_{i,j}$ which appears in the second term $\tilde W_{i,j}^p \left((W1)_i\right)$. The same cancellation appears for $\tilde W_{j,j} \tilde W_{i,j}^p$. Hence for $i\neq j$ the first two lines do not contain the diagonal terms of $W$.
\end{remark}
\begin{proof} The proof is a direct application of It\^o Formula. For \eqref{eq_tW_A}, let us first compute 
for $i',j' \in \cV$ and $i,j\in \cV$:
\begin{eqnarray*}
\left< d\overline M_{i',j'}, dM_{i,j} \right >&=& \sum_k \overline\sigma_{i',j'}^k X_k^2 \sigma_{i,j}^k  +\overline\sigma_{i,j}^k X_k^2 \sigma_{i',j'}^k =2\Re( \sum_k \overline\sigma_{i',j'}^k X_k^2 \sigma_{i,j}^k ).
\end{eqnarray*}
 Note that if $i,j=i',j'$, then 
\begin{eqnarray*}
\left< d\overline M_{i,j}, dM_{i,j} \right >&=& 2W_{i,j} (W X^2 W)_{i,j}  dt,
\end{eqnarray*}
Besides we have with the Ward identity~\eqref{Ward}
\begin{eqnarray*}
\left< dM_{i,j}, dX_j\right >&=& -{1\over \sqrt{2}} W_{i,j} X_j^2 ( \i G_V(j,j)-\i\overline{G_V(j,j)}) 
\\
&=&\sqrt{2}  W_{i,j} X_j^2 \Im (G_V(j,j))
\\
&=&  -W_{i,j} X_j^2 (WX)_j.
\end{eqnarray*}
since $\left< dV_k,dX_k\right>=  {\i\over \sqrt{2}} X_k^2 dt.$
Similarly,
$$
\left< dM_{i,j}, dX_i\right >=  -W_{i,j} X_i^2 (WX)_i.
$$
Hence it gives \eqref{eq_tW_A} from \eqref{eq_W_A}, 
where the second term comes from the square bracket between $dX_i$ and $dW_{i,j}$ and $dX_j$ and $dW_{i,j}$.

Let us prove \eqref{eq_tW_A_p}. The square brackets of the martingale terms are given by
\beq
&&\left< \tilde W_{i',j'}(dB_{i'}(t)+dB_{j'}(t)), X_{i}X_j \sum_{k} \left( \sigma^{ij}_k dV_k(t)+\overline{\sigma^{ij}_k dV_k(t)}\right) \right>
\\
&=&
- \sqrt{2}  \left( (\tW 1)_{i'}  X_i X_j \hbox{Im}( \sigma^{ij}_{i'}) +  (\tW 1)_{j'}  X_i X_j \hbox{Im}( \sigma^{ij}_{j'})\right)
\eeq
and since $\left< dV_k, dV_{k'}\right>=0 $ and  $\left< dV_k, d\overline {V_{k'}} \right>=\delta_{k=k'} $, we have,
\beq
&&\left< X_{i'}X_{j'} \sum_{k} \left( \sigma^{i'j'}_k dV_k(t)+\overline{\sigma^{i'j'}_k dV_k(t)}\right), X_{i}X_j \sum_{k} \left( \sigma^{ij}_k dV_k(t)+\overline{\sigma^{ij}_k dV_k(t)}\right) \right>
\\
&=&
\sum_{k} X_i X_j \overline{\sigma^{ij}_k} X_{i'} X_{j'} \sigma^{i'j'}_k + X_i X_j {\sigma^{ij}_k} X_{i'} X_{j'} \overline{\sigma^{i'j'}_k}
\eeq

Since  $\hbox{Im}( \sigma^{ij}_{i})=- \vert G_V(i,j)\vert^2 \Im(G_V(i,i))={1\over \sqrt{2}} W_{i,j} (WX)_i$, and $\vert \sigma^{ij}_k\vert^2=(W_{i,k}W_{k,j})W_{i,j}$, 
we deduce \eqref{eq_tW_A_p} from the previous identity for $i',j'=i,j$.
\end{proof}

\subsection{Case of a boundary set}\label{sec_SDE_A_bord}
Since $X_i(t)\to 0$ when $t\to \infty$, we see that the martingale $(W(t)X(t))_{i_0}$ or the variable $\tW_{i,j}(t)$ tend to $0$, hence we cannot recover interesting quantities from the limiting values. To solve this problem 
it is convenient to introduce a boundary set which leads to non trivial limiting values. 

Consider a subset $\cV'\subset \cV$ and we freeze the potential $V(t)$ on $\cV\setminus\cV'$. More precisely, we fix
\begin{eqnarray}\label{boundary_A}
\forall i\in \cV\setminus \cV', \;\;\;
\begin{cases}
X_i(t)=1, \forall t\ge 0,
\\
Y_i(t)=0, \forall t\ge 0.
\end{cases}
\end{eqnarray}
so that $V(t)={\i\over\sqrt{2}}$ on $ \cV\setminus \cV'$.
It is easy to see that for the boundary case, the equation \eqref{eq_W_A} is now
\begin{eqnarray}\label{eq_W_A_'}
\begin{cases}
dW_{i,j}(t)  = ( \sum_{k\in \cV'} W_{i,k} X_k^2 W_{k,j} ) dt+ dM_{i,j},
\\
dX_i(t)=\indic_{i\in \cV'} X_i(t) dB_i(t).
\end{cases}
\end{eqnarray}
on where $dM_{i,j} $ has the same expression as in \eqref{dM} where $dV_k=0$ for all $k\in \cV\setminus\cV'$.

Besides, looking at the proof of Proposition~\ref{prop_eq_tW_A}, and considering that $dV_k=0$ for all $k\in \cV\setminus \cV'$, we get that $\tilde W(t)$ satisfies now the equation, for all $i,j\in \cV$, 
\begin{eqnarray}\label{eq_tW_A_p_'}
&&d\tW_{i,j}^p
\\
\nonumber 
&=& p \tW_{i,j}^{p-1} (\tilde W\cdot\tW )_{i,j}dt  - p  \tilde W_{i,j}^{p} \left(\indic_{i\in \cV'} (\tW\cdot 1)_i+ \indic_{j\in \cV'} (\tW\cdot 1)_j\right)dt 
\\
\nonumber &&+ {p(p-1)\over 2} \tilde W_{i,j}^p(\indic_{i\in \cV'} + \indic_{j\in \cV'})+ p\tilde W_{i,j}^p(\indic_{i\in \cV'} dB_i(t)+ \indic_{j\in \cV'} dB_j(t))+p^2\delta_{i=j}\delta_{i\in \cV'} dW^p_{ii} dt
\\
\nonumber 
&&+ p \tW_{i,j}^{p-1} X_iX_j dM_{i,j}
\end{eqnarray}
where we define
\begin{equation}\label{WdotW}
(\tW\cdot \tW)_{i,j}=\sum_{k\in \cV'} \tW_{i,k}\tW_{k,j},
\end{equation}

\section{The Brownian interpolation of the $\H^{2|2}$ model}\label{s_SDE_S}
\subsection{Definition and properties of the Schr\"odinger operator associated with the $\H^{2|2}$ model on a finite graph.}\label{sec_def_beta}
Let $\cV$ be a finite set and
$$
W:=(W_{i,j})_{i,j\in \cV},
$$
be a symmetric matrix with non-negative coefficients. We intentionally use the same letter $W$ as in the previous section: as we will see later, $W$ will be the initial value $W(0)=W$ of a stochastic process $W(t)$ which has several similarities with the one defined in Section~\ref{s_SDE_A} for the Anderson model. We will write later $W^A(t)$ (for Anderson) and $W^S(t)$ (for SuSy) when there is an ambiguity.

We denote by \(\ggg=(\cV,\cE)\) the non-oriented graph induced by $W$, where $(i,j)$ is an edge if and only if $W_{i,j}>0$. However, the comparison  with the Anderson model will involve the complete graph on $\cV$.

For \(\beta=(\beta_x)_{x\in \cV} \in\R^\cV\), we define \(H_\beta=\beta-W\), where  \(\beta\) denotes here abusively the diagonal matrix with diagonal coefficients \((\beta_i)_{i\in \cV}\); in particular, \(H_\beta\) is a \(\cV\times \cV\) matrix. We write $H_\beta>0$ to mean that $H_\beta$ is positive definite and we denote
$$
G_\beta:=(H_\beta)^{-1},
$$
when $H_\beta$ is invertible.
We recall the following definition.
\begin{theoA}[Theorem 4 in \cite{STZ17}, Theorem 2.2 in \cite{Letac-Wesolowski}, Lemma~4 of \cite{SZ19}]\label{prop:nuwtheta}
For all \(\theta\in(\R_+^*)^\cV\) and \(\eta\in(\R_+)^\cV\), the measure \(\nu_\cV^{W,\theta,\eta}\) defined by
\[
\nu_\cV^{W,\theta,\eta}(d\beta)=\ind_{H_\beta>0}
\l(\frac{2}{\pi} \right)^{|\cV|/2}
\exp\l(-\frac{1}{2}\langle\theta,H_\beta\theta\rangle-\frac{1}{2}\langle\eta,G_\beta \eta\rangle+\langle\eta,\theta\rangle \r)
\frac{\prod_{i\in \cV}\theta_i}{\sqrt{|H_\beta|}}d\beta
\]
is a probability distribution. Moreover, for all \(i\in \cV\), the random variable \(\frac{1}{2\beta_i-W_{i,i}}\) has Inverse Gaussian distribution with parameter \((\frac{\theta_i}{\eta_i+\sum_{j\neq i}W_{i,j}\theta_j},\theta_i^2)\). 

\end{theoA}
The following property is elementary but important to understand the transformation we do below in Section~\ref{sec_Gauge_S}.
\begin{proposition}[Gauge transform]\label{prop_Gauge}
If $\beta$ is distributed according to $\nu_\cV^{W,\theta,\eta}$ then $\theta^2\beta$ has law $\nu_\cV^{\theta W\theta,1,\theta\eta}$ where we write $(\theta W\theta)_{i,j}=\theta_i W_{i,j}\theta_j$ and $(\theta\eta)_x=\theta_x\eta_x$.

We simply denote by $\nu_\cV^{W,\eta}$ the law $\nu_\cV^{W,1,\eta}$ when $\theta_x=1$ for all $x\in \cV$. When $\theta=1$ and $\eta=0$, we simply write $\nu_\cV^W$.
\end{proposition}
The potential $\beta$ was introduced in \cite{STZ17} without the boundary term $\eta$ and generalized later by Letac  (published later in \cite{Letac-Wesolowski}). In \cite{STZ17} it was shown to be related the first marginal of the supersymmetric hyperbolic sigma model $\H^{2|2}$ in horospherical coordinates. The $\H^{2|2}$-model was previously shown to be equivalent to a reinforced process, the VRJP (\cite{ST15}). It was in particular investigated in the case of the nearest neighbors lattice $\Z^d$ with constant conductances. 

In Section~\ref{sec_about_beta} we give extra properties of the potential $\beta$ and of the associated Green function $G_\beta$.
In particular, we introduce several key tools and recall the main results about the localization and delocalization results of the Green function $G_\beta$.

\subsection{The Stochastic Differential Equation associated with the $\beta$-potential}
We first recall some notations of the introduction. For $T=(T_i)_{i\in \cV}\in (\R_+)^V$ a vector of positive reals (which should be understood as times indexed by the vertices) we denote by $K^W_T$ or simply $K_T$ the matrix
$$
K_T=K_T^W=\Id-TW,
$$
where in the previous expression $T$ denotes the operator of multiplication by $(T_i)$.

The theorem below gives a representation of the distribution $\nu_\cV^{W,\eta}$ in terms of the final quadratic variation of a stochastic process. 
\begin{theorem}[Theorem~2.1 and Theorem~2.2 \cite{GRSZ24}]\label{thm_SDE}
(i) The following SDE admits a unique strong solution $X(t)=(X_i(t))_{i\in \cV,t\in \R_+}$ with initial condition $X_i(0)=1$,
\begin{eqnarray}\label{SDE_X}
\begin{cases}
dX_i(t)=X_i(t)dB_i(t)-X_i(t)^2\l(W(K^W_{T(t)})^{-1} X(t)\r )_i dt,
\\
T_i(t)=\int_0^t X_i(s)^2 ds=\int_0^t \l<dX_i(s),dX_i(s)\r>.
\end{cases}
\end{eqnarray}
More precisely, $T(t)$ remains such that $K^W_{T(t)}>0$, and $X_i(t)\ge 0$ for all $t$ and $i\in \cV$ and $\lim_{t\to \infty} X_i(t)=0$.

(ii) For all $i\in \cV$, the following limit exists in $\R_+$ and
$$
T_i:=\lim_{t\to\infty} T_i(t),
$$
is such that $1/T=(1/T_i)_{i\in \cV}$ has distribution $\nu_\cV^{W}$.
\end{theorem}
\begin{remark} Note that if we start the processus with the initial condition $X_i(0)=\theta_i>0$, then $1/T$ has law $\nu_\cV^{W,\theta}$. We explain below how to obtain the boundary term $\eta$.
\end{remark}
Remark now that 
$$
{\partial\over\partial t} (K^W_{T(t)})^{-1}= (K^W_{T(t)})^{-1} X(t)^2 W (K^W_{T(t)})^{-1}.
$$
Hence if we set
$$
W(t):= W(K^W_{T(t)})^{-1},
$$
we see that $(W(t),X(t))$ is a diffusion process satisfying the following SDE
\begin{eqnarray}\label{SDE_W_S}
\begin{cases}
dW(t)_{i,j}=(W(t) X(t)^2 W(t))_{i,j} dt,
\\
dX_i(t)=X_i(t) dB_i(t)-X_i(t)^2(W(t) X(t))_i dt.
\end{cases}
\end{eqnarray}

\subsection{The Gauge transform $\tW_{x,y}$ for the $\H^{2|2}$ model}\label{sec_Gauge_S}

As in Section~\ref{sec_Gauge_A} we set, for all $i,j\in \cV$,
$$
\tilde W_{i,j}(t)=X_i(t) W_{i,j}(t) X_j(t).
$$
This transformation corresponds to the normalization of the parameter $\theta$, see Proposition~\ref{prop_Gauge}.
 
Using \eqref{SDE_W_S},
we have by It\^o formula the following result.
\begin{proposition}\label{prop_tW_S}
The process $(\tilde W_{i,j}(t))_{i,j \in \cV}$ is a diffusion process solution of the following stochastic differential equation: for all $i,j\in \cV$,
\begin{eqnarray}\label{eq_tW_S}
&&
d\tilde W_{i,j} =
(\tilde W^2)_{i,j}dt  -\tilde W_{i,j} \left((\tW\cdot 1)_i+(\tW\cdot 1)_j\right)dt + \tilde W_{i,j}(dB_i(t)+dB_j(t))+\delta_{i=j} W_{ii} dt
\end{eqnarray}
More generally, for $p>0$, we have
\begin{eqnarray}\label{eq_tW_S_p}
d\tW_{i,j}^p&=& p \tW_{i,j}^{p-1} (\tilde W^2 )_{i,j}dt  - p  \tilde W_{i,j}^{p} \left((\tW\cdot 1)_i+(\tW\cdot 1)_j\right)dt 
\\
\nonumber &&+ p(p-1)\tilde W_{i,j}^p+ p\tilde W_{i,j}^p(dB_i(t)+dB_j(t))+ p^2 \delta_{i=j} W_{ii}^p dt
\end{eqnarray}
\end{proposition}
\begin{remark} 
Note that as in remark~\ref{rk_diag} the diagonal terms in the right-hand side of \eqref{eq_tW_S_p} cancel when $i\neq j$. 
Hence changing the definition of $\tilde W_{i,j}$ by putting $0$ on the diagonal, the identity above remains valid, and $(\tilde W_{i,j})_{i\neq j}$ is itself a diffusion process.
\end{remark}
\subsection{Boundary term}\label{sec_SDE_S_bord}
As in the Anderson case, since $X_i(t)\to 0$ when $t\to \infty$, we see that the martingale $(W(t)X(t))_{i_0}$ and $\tW_{i,j}(t)$ tends to $0$. In order to keep some information in the limit we introduce a model with a boundary subset.

A non trivial boundary term $\eta$ in $\nu_\cV^{W,\eta}$ appears naturally either when one adds a drift term as in \cite{GRSZ24} or by freezing some of the Brownian motions, which is more natural in our context. More precisely, we consider a subset $\cV'\subset \cV$ and suppose that for $i\in \cV\setminus\cV'$ we freeze the process $X_i(t)$, i.e., we have 
\begin{eqnarray}\label{X_bord}
\forall i\in \cV\setminus\cV',\;\;\; \begin{cases}
X_i(t)=1, \;\;\; \forall t\ge 0, 
\\
T_i(t)=0, \;\;\; \forall t\ge 0.
\end{cases}
\end{eqnarray}
We define $K_{T(t)}^W$ as before, which means that  $K_{T(t)}^W$ has a 1 on the diagonal and $0$ outside on any line $i\in \cV\setminus\cV'$.
Then the theorem~\ref{SDE_X} has the counterpart below
\begin{theorem}[Theorem~2.1 and Theorem~2.2 \cite{GRSZ24}]\label{thm_SDE_eta}
(i) The SDE \eqref{SDE_X} with boundary values \eqref{X_bord} on $\cV\setminus\cV'$ has a unique solution $(X_i(t))_{i\in \cV'}$.

(ii) For all $i\in \cV'$, the following limit exists in $\R_+$ and
$$
T_i:=\lim_{t\to\infty} T_i(t),
$$
is such that $1/T=(1/T_i)_{i\in \cV'}$ has distribution $\nu_{\cV'}^{W_{ \cV', \cV'},\eta}$ where $W_{ \cV', \cV'}$ is the restriction of the matrix $W$ to $\cV'$ and $\eta=W_{ \cV',\cV\setminus\cV'} \indic_{\cV\setminus\cV'}$, i.e. $\eta_i=\sum_{j\in \cV\setminus\cV'} W_{i,j}$ for $i\in \cV'$.
\end{theorem}
\begin{proof} This is direct consequence of Theorem~2.1 and Theorem~2.2 \cite{GRSZ24} since if we freeze the process $X_i(t)$ for $i\in \cV\setminus \cV'$, then by writing $K_T$ in terms of its restriction to $\cV'\times \cV'$, the identity on $\cV\setminus\cV'\times \cV\setminus\cV'$ and the cross term $(-T(t) W)_{\cV'\times \cV\setminus\cV'}$, the SDE \eqref{SDE_X} with boundary values \eqref{X_bord} is equivalent to the SDE of Theorem~2.1 \cite{GRSZ24} with drift $\eta$.
\end{proof}

With the definition above, we see that ${\partial\over\partial t} T_i(t)= \indic_{i\in \cV'} X_i(t)^2$. Hence, we easily see that $(W(t),X(t))$ satisfies for all $i,j\in \cV$,
\begin{eqnarray}\label{SDE_W_S_bord}
\begin{cases}
dW(t)_{i,j}=\sum_{k\in \cV'} W_{i,k} (t) X_k(t)^2 W_{k,j}(t)) dt,
\\
dX_i(t)=\indic_{i\in  \cV'} (X_i(t) dB_i(t)-X_i(t)^2(W(t) X(t))_i dt).
\end{cases}
\end{eqnarray}
It easily implies that $\tilde W(t)$ satisfies now the equation, for all $i,j\in \cV$, 
\begin{eqnarray}\label{eq_tW_S_p_'}
&&
d\tW_{i,j}^p
\\
\nonumber 
&=& p \tW_{i,j}^{p-1} (\tilde W\cdot\tW )_{i,j}dt  - p  \tilde W_{i,j}^{p} \left(\indic_{i\in \cV'} (\tW\cdot 1)_i+ \indic_{j\in \cV'} (\tW\cdot 1)_j\right)dt 
\\
\nonumber &&+ p(p-1)\tilde W_{i,j}^p(\indic_{i\in \cV'} + \indic_{j\in \cV'})+ p\tilde W_{i,j}^p(\indic_{i\in \cV'} dB_i(t)+ \indic_{j\in \cV'} dB_j(t))+p^2\delta_{i=j}\delta_{i\in \cV'} dW^p_{ii} dt
\end{eqnarray}
where  $\tW\cdot \tW$ is defined in \eqref{WdotW}.
\subsection{Comparison of the equations for the Anderson model with Cauchy potential and the $\H^{2|2}$.}\label{sec_comparison}
As we can see in \eqref{eq_tW_A}, \eqref{eq_tW_A_p} and \eqref{eq_tW_S}, \eqref{eq_tW_S_p} (and also for the corresponding equations with boundary values), there are strong formal similarities between the two sets of equations. More precisely, for $p=1$, there is only an extra martingale term $X_iX_j dM_{i,j}$ in \eqref{eq_tW_A} compared to \eqref{eq_tW_S}. However, this extra martingale term has a non trivial influence, in particular it is responsable for the difference in the equations for $\tW_{i,j}^p$. Indeed, for the equation for a general $p$, we see that the first line has a factor $p^2$ in \eqref{eq_tW_A_p} and a factor $p$ in \eqref{eq_tW_S_p}, while the second lines are the same. This means that the quadratic variation due to the extra term martingale $X_iX_j dM_{i,j}$ is of the order of the terms in the first line.

{\bf Initial condition.} The natural question is whether we can compare the final values of the solutions of the two equations. We denote by $\tW^A(t)$ and $\tW^S(t)$ the solutions of respectively the Anderson equation \eqref{eq_tW_A} and the $\H^{2|2}$ equation \eqref{eq_tW_S}. The value of $\tW^A(0)$ is fixed by the model and is equal to $\tW_{i,j}^A(0)=\vert (({\i\over \sqrt{2}} \Id +{1\over \lambda}( -\Delta -E))^{-1})_{i,j} \vert^2$.

We remark that in \eqref{eq_tW_A_p} and \eqref{eq_tW_S_p}, the first line is $p+1$-homogeneous in $\tW$ and the second line is $p$-homogeneous. It means that defining $\hat W^S:= {1\over p} \tW^S$,  $(\hat W^S)^p$ would exactly satisfies the equation given by the first two lines of the equation for $(\tW^A)^p$ in \eqref{eq_tW_A_p} (hence with the martingale term $dM_{i,j}$ removed). Hence, considering the equations for $p$ would give a natural choice of initial value $\tW^S(0)= p\tW^A(0)$. But if we want to compare the localization/delocalization of the two models it is not clear which $p$ we should take.

In the next section we analyse the phase diagram of the $\H^{2|2}$ model with conductances given by $r \tW_{i,j}^A(0)$ on $\Z^d$, $r>0$, and show that whatever $r$, the diagram displays a localized and delocalized phase similar to the ones proved or expected for the Anderson model with Cauchy potential.

In the last section~\ref{sec_discussion}, we discuss further the relations between the two equations and propose a conjecture.

\section{About the phase diagram of the associated $\mathbb{H}^{2 | 2}$.}\label{sec_delocalization}

\subsection{Complementary material about the $\beta$ potential}\label{sec_about_beta}
We consider a finite set $\cV$ and a symmetric matrix $W=(W_{i,j})_{i,j\in \cV}$ with non negative coefficients as in Section~\ref{sec_def_beta}. Recall the definition of $\nu_\cV^{W,\eta}$ and $\nu_\cV^W$ in Theorem~\ref{prop:nuwtheta} and Proposition~\ref{prop_Gauge}. The following proposition will be useful.

\begin{proposition}[Theorem~3 of \cite{STZ17}]\label{law_gamma}
If $\beta$ is distributed under $\nu_\cV^W$, then for all $i\in \cV$, $G_\beta (i,i)$ has the same law as ${1\over 2 \gamma}$ where $\gamma$ is a Gamma random variable with parameters $(\demi, 1)$. Besides, $G_\beta (i,i)$ is independent of $(\beta_j)_{j\neq i}$.

\end{proposition}

For all $i_0\in \cV$, the function
$$
u_{i_0}(j):=\ln\left({G_\beta(i_0,j)\over G_\beta(i_0,i_0)}\right),
$$
has the law of the first marginal of the $\H^{2|2}$ supersymmetric sigma field in horospherical coordinates defined in \cite{DSZ10} rooted at $i_0$, see Theorem~3 of \cite{STZ17} for a precise statement. In \cite{DSZ10,DS10}, it is shown, on large boxes of $\Z^d$ with constant nearest neighbor conductances,  that the field $e^{u_0(j)}$ is exponentially decreasing for small conductances (strong disorder) and for $d\ge 3$, delocalized for large conductances  (weak disorder) in the sense that $\E\left( e^{m(u_0(i)-u_0(j))}\r)\le 2$ for some $m>0$ and all vertices $i,j$. We make this more precise in Section~\ref{sec_loc_deloc} and introduce the critical parameters.

The $\beta$ representation is in particular useful since it displays interesting properties with respect to restriction. The following Lemma was found independently in \cite[Lemma~5]{SZ19} and \cite[Propositions~4.2 and~4.3]{Letac-Wesolowski} and will be useful later.
\begin{lemA}
\label{lem:conditionalbeta} 
Assume $\cV$ is finite and let $U\subset \cV$. Under $\nu_\cV^{W,\eta}(d\beta)$, the following holds:
\begin{itemize}
  \item[(i)] The restriction $\beta_U$ is distributed according to $\nu^{W_{U,U},\hat\eta}_U$, where
  \begin{equation}
  \hat\eta:= \eta_U + W_{U,U^c} 1_{U^c}\,.
  \end{equation}
  \item[(ii)] Conditionally on $\beta_U$, its complementary part $\beta_{U^c}$ is distributed according to $\nu^{\check{W},\check{\eta}}_{U^c}$, where $\check{W}=(\check{W}_{i,j})_{i,j\in U^c}$ and $\check\eta\in(\R)^{U_c}$ are defined by
  \begin{equation}\label{eq:lem:conditionalbeta}
    \begin{split}
  \check{W} & := W_{U^c,U^c} + W_{U^c,U}((H_\beta)_{U,U})^{-1}W_{U,U^c} \\
   \text{and}\qquad \check\eta & := \eta_{U^c}+W_{U^c,U}((H_\beta)_{U,U})^{-1} \eta_U\,.
    \end{split}
  \end{equation}
\end{itemize}
\end{lemA}
The previous lemma will be useful to prove the property of block reduction of the $\beta$ potential.
A similar property was used in \cite{Disertori-Merkl-Rolles} to analyse a hierarchical version of the $\H^{2|2}$ model. However, we present a different approach using the random Schr\"odinger representation rather than supersymmetry.
\begin{lemma}\label{lem_block}
Let $\cV$ be finite and decomposed as
\begin{equation}\label{decomposition}
\cV=\sqcup_{x\in X} B_x,
\end{equation}
where $(B_x)_{x\in X}$ forms a partition in non-empty subsets of $\cV$. Let $W=(W_{i,j})_{i,j\in \cV, }$, $\eta=(\eta_i)_{i\in \cV}$ be some conductances and boundary terms such that the matrix $W$ and the vector $\eta$ are constant on blocks of the partition above \eqref{decomposition}. We define
 $\tilde W=(\tilde W_{x,y})_{x,y\in X, x\neq y}$, $\tilde \eta=(\tilde \eta_x)_{x\in X}$ such that
 \begin{eqnarray}  \label{W_tilde}
\tilde W_{x,y} = \vert B_x\vert  W_{i,j} {\vert B_y\vert }, \;\;\; \hbox{for all $x,y\in X$, $x\neq y$, $i\in B_x, \; j\in B_y$}
 \end{eqnarray}
 and
 \begin{eqnarray} \label{eta_tilde}
 \tilde \eta_{x}={\vert B_x\vert } \eta_{i}, \;\;\; \hbox{for all $x\in X$, $i\in B_x$}.
 \end{eqnarray}
 Define
 $$
 \hat G^{B_x}=((\beta-W)_{B_x\times B_x})^{-1},
 $$
 where $\beta$ is distributed according to $\nu_{\cV}^{W,\eta}$ and set
 \begin{eqnarray}
 \label{beta_tilde}
  \tilde \beta_x={ \vert B_x\vert^2 \over \left< \indic_{B_x} , \hat G^{B_x} \indic_{B_x}\right>}.
  \end{eqnarray}
 Then $\tilde \beta$ has law $\nu_{X}^{\tilde W,\tilde \eta}$. Besides, for $x_0\in X$,
 \begin{eqnarray}\label{block_psi}
  {1\over \vert B_{x_0}\vert} \left< \indic_{B_{x_0}} , G \eta\right>
 \end{eqnarray}
 has the law of $\tilde \psi_{x_0}:= \tilde G \tilde \eta$ under $\nu_{X}^{\tilde W,\tilde \eta}$.
  \end{lemma}
In order not to delay too much the main statement, the proof of this lemma will be given in Appendix~\ref{appendix_A}.

\subsection{Extension to infinite graphs and monotonicity principle}\label{sec_martingale}
From Lemma~\ref{lem:conditionalbeta} we can extend the definition of the $\beta$ potential to countable sets $\cV$. 
This was done in \cite[Section~2.2]{SZ19} in the case of local conductances, and generalized to non local conductances in \cite{BLPS25}.

Let $\cV$ be a countable set and $W=(W_{i,j})_{i,j \in \cV}$ be non-negative conductances such that $W_{i,j}=W_{j,i}$ and
\begin{equation}\label{finite_sum}
\sum_{j\in \cV} W_{i,j}<\infty, \;\;\; \forall i\in \cV.
\end{equation}
The following is an extension of \cite[Section~2.2]{SZ19} and was proved in \cite[Definition~C]{BLPS25}. 
\begin{propA}\label{prop_beta_infinite}
Let $\eta=(\eta_i)_{i\in \cV}\in \R_+^\cV$. There exists a unique probability distribution $\nu_\cV^{W,\eta}$ on $\R^\cV$ such that for all finite subset $U\subset \cV$, $\beta_{U}$ has distribution $\nu_{U}^{W_{U,U},\hat \eta_U}$, where
$$
\hat \eta_U=\eta_U + W_{U,U^c} 1_{U^c} \,.
$$
\end{propA}

Consider now $\eta=0$ and $\beta = (\beta_x)_{x\in \cV}$ distributed according to $\nu_\cV^{W}$. We define as before
$$
H_\beta=\beta-W\,,
$$
where the operator $W$ on $\R^\cV$ is defined by $Wf(x)=\sum_{y\in \cV} W_{x,y} f(y)$ if the function $f$ verifies $\sum_{y\in \cV} W_{x,y} \vert f(y)\vert<\infty$ for all $x$. 
Note that by Theorem~\ref{prop:nuwtheta} and Proposition~\ref{prop_beta_infinite}, for any finite subset $U\subset \cV$ we have $\left(H_\beta\right)_{U,U}>0$. Consider now an increasing sequence $(\cV_n)_{n\ge 0}$ of finite subsets such that $\cup_n \cV_n=\cV$, and define  
$$
\hat G^{(n)}_\beta= \left( \left(H_\beta\right)_{\cV_n,\cV_n}\right)^{-1} 
$$ 
the Green function of the operator $H_\gb$ restricted to $\cV_n$.
Since $\left(H_\beta\right)_{\cV_n,\cV_n}>0$, $\hat G^{(n)}_\beta$ has the following path expansion,
\begin{equation}\label{G_paths}
\hat G^{(n)}_\beta(x,y)= \sum_{\sigma \colon x\xrightarrow{\cV_n} y} \frac{ \prod_{i=0}^{\vert \sigma\vert -1} W_{\sigma_i,\sigma_{i+1}}}{\prod_{i=0}^{\vert \sigma\vert} \beta_{\sigma_i}}. 
\end{equation}
where the sum on $\sigma \colon x\xrightarrow{\cV_n} y$ runs on all finite paths $\sigma$ on the graph $\ggg$ from $x$ to $y$, which stay inside $\cV_n$.

We also define $(\psi_{n}(x))_{x\in \cV}$ as the solution of
\begin{equation}\label{eq:def:psin}
\begin{cases}
\psi_{n}(x)=1,& \;\;\; \hbox{if $x\notin \cV_n$}\,,
\\
H_\beta(\psi_{n})(x)=0,& \;\;\; \hbox{if $x\in \cV_n$}\,.
\end{cases}
\end{equation}
which is unique since $(H_\beta)_{\cV_n,\cV_n}>0$. Besides, if we denote $\hat \eta_{\cV_n}= W_{\cV_n, \cV\setminus\cV_n}$, i.e. $\forall x\in \cV_n$, $\hat \eta_{\cV_n}(x)=\sum_{y\in \cV_n^c} W_{x,y}$, we have
\begin{eqnarray}\label{Gbeta}
\psi_{n}(x)=\hat G^{(n)}_\beta\hat \eta_{\cV_n}.
\end{eqnarray}
Note that $\psi_n(x)$ is nothing else than the first marginal of the $\H^{2|2}$-model for a rooting at $\delta$ the external point of wired boundary condition outside of $\cV_n$. 
The interest of considering $\psi_{n}$ is in its martingale property and in that it characterizes the recurrence/transience of the associated reinforced process, the VRJP.

The following is proved in \cite[Theorem~1]{SZ19}.

\begin{theoA}\label{Th_martingale}
For all $x,y\in \cV$, the sequence of random variables \((\hat{G}^{(n)}(x,y))_{n\geq 0}\) is non decreasing and converges a.s.\ to 
\begin{eqnarray}\label{Ghat}
\hat G_\beta(x,y):= \lim_{n\to \infty} \hat G_\beta^{(n)}(x,y).
\end{eqnarray}
Moreover, $\nu_\cV^W$-almost surely, $0<\hat G_\beta(x,y)<\infty$ and the limit does not depend on the choice of the sequence of subsets $\cV_n$. 

Additionally, under the probability $\nu_\cV^W$, for all \(x\in \cV\), \((\psi_{n}(x))_{n\geq 0}\) is a positive \(\cF_{n}\)-martingale with $\cF_n=\sigma\{\beta_x, \; x\in \cV_n\}$. It converges a.s.\ to a
random variable \(\psi(x)\), such that $\psi(x)\ge 0$ a.s., and the limit does not depend on the choice of the
increasing sequence \((\cV_n)_{n\geq 0}\).
Besides, the quadratic variation of the vectorial martingale \((\psi_{n}(x))_{x\in \cV}\) is $\hat G^{(n)}(x,y)$ and
\(\psi_{n}(x)\) is bounded in~\(L^2\) if and only if \(\E_{\nu_\cV^W}(\hat G(x,x))<\infty\).

If the conductances are vertex transitive (e.g. if $\cV=\Z^d$ and $W_{i,j}$ depends only on $j-i$) then either $\psi=0$ a.s. in which case the  VRJP associated with conductances $W$ is recurrent, or $\psi>0$ a.s. in which case the VRJP associated with conductances $W$ is transient.
\end{theoA}
Note that by construction, the function $\psi$ satisfies $H_\beta \psi=0$ and that $\hat G_\beta$ is a Green function in the sense that $H_\beta \hat G_\beta=\Id$. 

The operator $\hat G_\beta$ is a Green function of $H_\beta$ on $\cV$ in the sense that it satisfies $H_\beta \hat G_\beta=\Id$ but it is not the natural Green function of the model since it does not describe the VRJP nor it is the limit of the Green functions on finite boxes without boundary term. The natural definition is, see \cite[Theorem~1~(iii)]{SZ19}
\begin{eqnarray}\label{Ggamma}
G_{\beta,\gamma}(x,y):=\hat G_\beta (x,y)+\frac{1}{2\gamma} \psi(x)\psi(y),\;\; x,y\in \cV
\end{eqnarray}
where $\gamma$ is a $\mathrm{Gamma}(\demi,1)$ random variable, independent of $\beta$. Then $G_{\beta,\gamma}$ gives a representation of the VRJP as mixture of Markov jump processes, and $G_{\beta,\gamma}$ is the limit of the Green functions with free boundary conditions and no boundary term $\eta$. Besides, $G_{\beta, \gamma}$ satisfies the Proposition~\ref{law_gamma}, which will be useful later.

The following result of Poudevigne \cite[Theorem~6]{Poudevigne19} is a key tool to analyse the Green function of the $\beta$-potential and the $\H^{2|2}$ model. We state it in the form we will need, see \cite[Theorem~E]{BLPS25}

\begin{theoA}\label{thm_monotonicity}
Let $\ggg=(\cV,\cE)$ be a countable graph, and consider $W^-$ and $W^+$ some conductances as in \eqref{finite_sum}, such that for all $x,y\in \cV$, $W^-_{x,y}\le W^+_{x,y}$. 
Let $U\subset \cV$ a finite subset and $\eta^\pm=W^\pm_{U,U_c}\indic_{U^c}$. Then, there is a coupling $(\beta^-,\beta^+)$ with $\beta^\pm \sim \nu_U^{W^\pm,\eta^\pm}$, such that, if $\psi^\pm=G_{\beta^\pm}^{W^\pm}\eta^\pm$ as in definition \eqref{Gbeta}, then for  all $x\in U$ which is connected to $U^c$ by the conductances $W^-$ (and hence $W^+$), 
$$
\E\left[ \psi^-(x) \mid \beta^+\right]=\psi^+(x),
$$
In particular, for any non-negative convex function~$f$, we have that 
\[
  \E\left[ f\big(\psi^-(x) \big) \mid \beta^+\right] \geq f\big( \psi^+(x) \big) \,.
\]
\end{theoA}

\subsection{Localization and delocalization for the nearest neighbor model}\label{sec_loc_deloc}
The martingale $\psi_n$ and the monotonicity Theorem~\ref{thm_monotonicity} allow to give a clear definition of the critical values of the model.
Consider now the case of the lattice $\Z^d$ with nearest neighbor conductances $W$. By monotonicity,  there exists a unique critical value $W_c(\Z^d)$ such that $\psi=0$ a.s. for $W<W_c(\Z^d)$ and $\psi>0$ a.s. for $W>W_c(\Z^d)$, see \cite{Poudevigne19} Theorem~2. The localization result of \cite{DS10} imply $W_c(\Z^d)>0$ for all $d\ge 1$ and that the $\H^{2|2}$ model is exponentially localized for $W$ small enough (which can be expressed equivalently in terms of exponential decay of the $\H^{2|2}$ field $e^{u_0}$, of $G_\beta(x,y)$ or of $\psi_n$). For $d=1$ and $d=2$, $W_c(\Z^d)=\infty$ (by \cite{Davis_Volkov,SZ19,Sabot21,Kozma_Peled21}). For $d\ge 3$, the important delocalization result of \cite{DSZ10} implies $W_c(\Z^d)\in (0,\infty)$. 

The existence of a delocalized phase on $\Z^d$ with nearest neighbor and $d\ge 3$, and the monotonicity will be useful in the proof of Theorem~\ref{thm_delocalization} below.

\subsection{Main theorem about the phase diagram of the $\H^{2|2}$ model with conductances given by $rW_{x,y}$ where $W_{x,y}=\l \vert \l({1\over \lambda} (-\Delta -E)+{\i\over \sqrt{2}} \r)^{-1}_{x,y}\r\vert^2$}

Fix a spectral value $E\in \R$, $\lambda>0$ the disorder parameter and a scale parameter $r>0$. On $\Z^d$ we consider the following infinite range weights: for all $x,y\in \Z^d$
\begin{eqnarray}\label{W_loc}
W_{x,y}:=\l \vert \l({1\over \lambda} (-\Delta -E)+{\i\over \sqrt{2}}\r)^{-1}_{x,y}\r\vert^2,
\end{eqnarray}
and we will also consider the scaled weights $(rW_{x,y})$.

The theorem below identifies a delocalized phase and a localized phase for the $\mathbb{H}^{2 | 2}$-model associated with weights $rW$ on $\Z^d$, $d\ge 3$. It relies heavily on the tools developed for the $\mathbb{H}^{2 | 2}$-model: localization, delocalization regimes for the nearest neighbor model on $\Z^d$ (\cite{DSZ10,DS10}), random Schr\"odinger operator representations and conditioning, restriction properties (\cite{STZ17,SZ19}), monotonicity principle (\cite{Poudevigne19}).

For technical reasons, the proof of delocalization is done for $E\in \sigma(-\Delta)$ minus a finite set of points where there are degeneracies on the Fermi curve or symmetries. We think the proof could be extended to $E\in ]0,4d[$.
\begin{theorem}\label{thm_delocalization}
(i) {\bf (Delocalization.)} Let $d\ge 3$. For all $E\in ]0,4d[ \setminus \sss_0$ where $\sss_0=\{4n, \;0\le n \le d\}\cup\{2d\}$, for all $r>0$ 
there exists $\lambda_0(E,r)>0$, such that for all $\lambda< \lambda_0(E,r)$, the $\mathbb{H}^{2 | 2}$-model with weights $rW$ is delocalized in the sense that $\psi_\infty>0$ a.s. with definitions of Section~\ref{sec_loc_deloc}. Besides, for $\lambda< \lambda_0(E,r)$ we can guaranty a uniform lower bound of $\E((\psi_\infty)^p)>c>0$ for $0<p<1$, locally in $E$.

(ii) {\bf (Localization.)} Let $r>0$. Under one of the following conditions:
\begin{itemize}
\item
for $E\not\in \sigma(-\Delta)$, and $\lambda <\lambda_1(E,r)$ for some $\lambda_1(E,r)>0$,
\item
for $\operatorname{dist}_\C (E+{\i\over \sqrt{2}}\lambda,\sigma(-\Delta) )>\Lambda(r)$, for some $\Lambda(r)>0$,
\end{itemize}
the $\H^{2|2}$ model with weights $rW$ is localized in the sense that $\psi_\infty=0$ a.s., with exponential decay of the Green function $\hat G_\beta$, more precisely, there exists $C>0$ and $\mu>0$ such that
$$
\E\left( (\hat G_\beta(0,x))^{1\over 4}\r )\le C e^{-\mu \| x \| },
$$
for $\hat G_\beta$ defined on $\Z^d\times \Z^d$ by \eqref{Ghat}. In terms of the martingale $\psi_n=\psi_n(0)$,
$$
\E\l ( (\psi_n)^{1\over 4}\r )\le C e^{-\mu n }.
$$
 \end{theorem}
 \begin{remark} The exponent ${1\over 4}$  in~{\it (ii)} is not important, it could be replaced by $0<s<\demi$ for $\hat G_\beta(0,x)$ (we need $s<\demi$ since, by Proposition~\ref{prop_Gauge}, $\hat G_\beta(0,0)^s$ is integrable only up to $\demi$) and by $0<s<1$ for $\psi_n$. Note also that since $\psi_\infty=0$, we have $\hat G_\beta=G_{\beta,\gamma}$ defined in \eqref{Ggamma}. 
 \end{remark}
\begin{figure}[htbp]
  \centering
  \includegraphics[width=\textwidth]{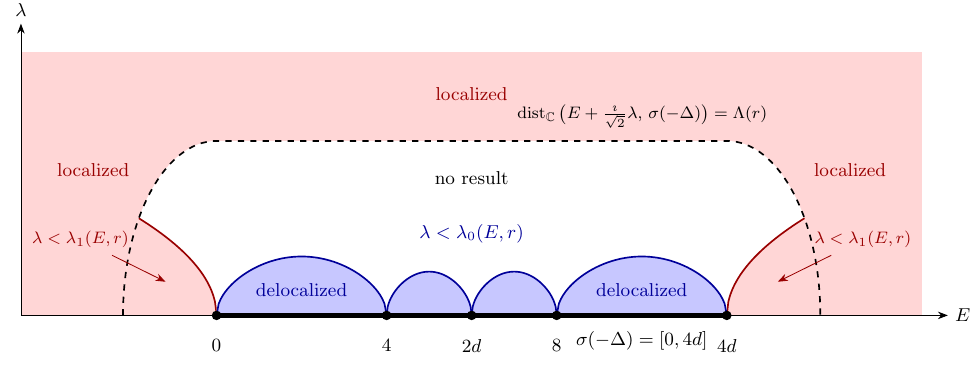}
  \caption{Phase diagram of Theorem~\ref{thm_delocalization} in the $(E,\lambda)$ plane ($d=3$).
}\label{figure1}
  \label{fig:phases}
\end{figure}
 \begin{proof}[Proof of Theorem~\ref{thm_delocalization}~(i)]
 For an integer $L>0$, and $x\in \Z^d$, we set
  $$C^{(L)}_x=(L x+]-L/2, L/2 ]^d)\cap \Z^d,
  $$
  so that $(C^{(L)}_x)_{x\in \Z^d}$ forms a partition of $\Z^d$ in boxes of size $L$.
We first prove the following.
 \begin{lemma}\label{lem_Green}
 Let $\lambda >0$ and $ E\in \sigma(-\Delta)\setminus\mathcal{S}_0$.
There exists $C>0$, $\epsilon >0$ such that  for $\lambda$ small enough, there exist an integer $L\ge {1\over \lambda}$ and subsets $B_x\subset C^{(L)}_x$ for all $x\in \Z^d$ such that, $0\in B_0$, 
 \begin{eqnarray}
 \label{card}
 \forall x\in \Z^d,\;\;\; &\operatorname{card}(B_x)=\operatorname{card}(B_0)\ge {\epsilon {\lambda}^{-d}}
 \\
 \label{lower_bound}
  &W_{i,j}\ge {C \lambda^{{d}+1}}, \;\;\; \forall i\in B_x, \; j\in B_y, \; \forall x\sim y,
\end{eqnarray}
  where $x\sim y$ means that $x$ and $y$ are nearest neighbors in the lattice $\Z^d$.
  \end{lemma} 
Let us first explain how we can derive Theorem~\ref{thm_delocalization}~(i) from that Lemma. Fix $\lambda $ and $L$ as above. The idea is first to use monotonicity to compare the random Schr\"odinger operator with distribution $\nu_{\Z^d}^{rW}$ with the one where $W_{i,j}$ is replaced by its lower bound \eqref{lower_bound} given in Lemma~\ref{lem_Green} when $i$ and $j$ are in adjacent boxes and $0$ otherwise, then to use Lemma~\ref{lem_block} and the results about delocalization for the $\H^{2|2}$ model. Hence we define for $i,j\in \Z^d$,
 $$
 \overline W_{i,j}=
 \begin{cases}
 {C \lambda^{d+1}},& \hbox{ if there exist $x\sim y$ in $\Z^d$ such that $i\in B_x$ and $j\in B_y$}
 \\
 0,& \hbox{ otherwise.}
 \end{cases}
 $$
Hence we have $W_{i,j}\ge \overline W_{i,j}$.

Consider the sequence of subsets $\cV_n=]-nL -L/2, nL+ L/2]^d \; \cap \Z^d$ and consider $\psi_n=\psi_n(0)$ the associated martingale as defined in Section~\ref{sec_martingale}. By Theorem~\ref{thm_monotonicity}, and by Lemma~\ref{lem_Green}, for $0<p < 1$,
$$
\E(\psi_n^p)\ge \overline \E(\psi_n^p),
$$
where on the right-hand side the expectation is with respect to the law $\nu_{\Z^d}^{\overline W}$ on $\Z^d$.
Besides, we can apply Lemma~\ref{lem_block} to the right hand side since $\overline W$ is constant between any two blocks $C_x^{(L)}$, hence by restriction to $\cV_n$, the boundary terms are also constant on blocks. 
More precisely, we can define $\overline{\overline W}_{x,y}$, for $x,y\in \Z^d$, by
$$
\overline{\overline W}_{x,y}=\indic_{x\sim y}\vert B_x\vert \vert B_y\vert  {C \lambda^{d+1}} 
\ge \epsilon^2 C {\lambda^{-2d} \lambda^{{d+1}}}\ge \epsilon^2 C \lambda^{-d+1}.
$$
Hence, we see that for all $r>0$ and $0<p<1$, for $\lambda $ small enough we have $r \overline{\overline W}_{x,y}> W_c^{(p)}(\Z^d)$ for $d\ge 3$. By Lemma~\ref{lem_block}, \eqref{block_psi}, applied to the conductances $\overline{\overline W}$
$$
\E(\psi_n(0)^p)\ge \overline E(\psi_n(0)^p)\ge{\vert B_0\vert^{p-1} }  \overline E\l (\l({\sum_{i\in B_0}  \psi_n(i)\over \vert B_0\vert}\r )^p\r) = {\vert B_0\vert^{p-1} }  {\overline {\overline E} }\l( (\psi_n(0)\r)^p),
$$
where in the last term the expectation is with respect to the law $\nu^{r\overline{\overline W}}$ restricted to ${\overline{\overline \cV}}_n= \{-n, \ldots, n\}^d$. Since $r \overline{\overline W}_{x,y}> W_c^{(p)}(\Z^d)$, it implies that $\psi_n(0)$ is bounded in $L^p$ for the conductances $r W$ given in \eqref{W_loc}.
Besides, since the controls can be taken uniform in a neighborhood $I$ of $E$, we can take a uniform lower bound $r \overline{\overline W}_{x,y}$ which implies a uniform lower bound on $E'\in I$ for  $\E((\psi_\infty)^p)$ for $0<p<1$.
  \end{proof}
 
\begin{proof}[Proof of lemma~\ref{lem_Green}]
Fix an energy $E\in ]0,4d[ \setminus \{4n, \;0\le n \le d\}$. We will explain later why we besides assume $E\neq 2d$ in Theorem~\ref{thm_delocalization}~(i).
For $\lambda>0$ a positive real we set
$$
g_{{\i\over \sqrt{2}}\lambda -E}=\left( -\Delta-E+{\i\over \sqrt{2}} \lambda\r)^{-1}={\lambda^{-1}} G_{{\i\over \sqrt{2}},E,\lambda},
$$
where $G_{V,E,\lambda}$ is defined in \eqref{GV},
so that we have
$$
W_{i,j}=\lambda^2 \l\vert g_{{\i\over \sqrt{2}\lambda}-E }(i,j) \r\vert^2.
$$
Written in Fourier, we have
\begin{eqnarray}\label{Fourier}
    g_{{\i\over \sqrt{2}}\lambda -E}(i, j) = \frac{1}{(2\pi)^d} \int_{\mathbb{T}^d} \frac{e^{i k (j-i)} }{\varepsilon(k) + {\i\over \sqrt{2}}\lambda -E} \, \mathrm{d}k,
\end{eqnarray}
where $\varepsilon(k) =2\sum_{\ell=1}^d (1-\cos(k_\ell))$ and $\mathbb{T}^d=(\R/ 2\pi \Z)^d=[0,2\pi[^d$.

Define the Fermi manifold as the set where the real part of the denominator vanishes, i.e.
$$
\Sigma_E=\{ k\in \mathbb{T}^d, \;\;\; \varepsilon(k)=E\}.
$$
The choice $E\not\in\{4n, \;0\le n \le d\}$ implies that $\vert \nabla \varepsilon(k)\vert \neq 0$ for all $k\in \Sigma_E$. Indeed, the gradient is $\nabla \varepsilon(k)=2(-\sin(k\cdot e_\ell))_{\ell=1,\cdots, d}$. If $\nabla \varepsilon(k)=0$, it implies that $\cos(k_\ell)=\pm 1$ for all $\ell=1, \cdots, d$, hence $E\in \sss_0$. 

We denote by $\S^d$ the Euclidian unit ball, and for $x\in \R^d\setminus\{0\}$, we set $\hat x={x\over \| x\| }\in \S^d$. For $\hat x\in \S^d$, by Lemma~3 of \cite{Vainberg2001}, there exists at most $2^d d$ points on $\Sigma_E$, that we denote $(k_j(\hat x))_{j=1, \ldots, m(\hat x)}$, such that $\nabla \varepsilon(k_j(\hat x))$ is parallel to and in the direction of $\hat x$. The vector $\hat x$ is non singular if the total curvature $K(k)$ of $\Sigma_E$ does not vanish on $k_j(\hat x)$ (see \cite{Vainberg2001}, proof of Lemma~2 in the appendix for a formula for $K(k)$). The set of non singular vectors is an open subset of $\S^d$.

The following lemma is an adaptation to $\lambda >0$ of the main result of \cite{Vainberg2001}.
\begin{lemma}\label{lem_Green_estimate}
Let $K>1$ be real. For $\lambda>0$ and $x\in \Z^d$ such that $\hat x$ is non singular and ${1\over K \lambda}\le \| x\| \le {K\over \lambda}$, we have the following estimate
$$
g_{{\i\over \sqrt{2}}\lambda -E}(0,x)=\sum_{j=1}^{m(\hat x)} 
{a_j(\hat x) }e^{\i k_j(\hat x)\cdot x}
 {e^{-\|x\| \lambda  \over  \sqrt{2} \vert \nabla \varepsilon(k_j(\hat x))\vert}} \| x\| ^{-{d-1\over 2}}+ O(\lambda^{{d+1\over 2}}),
$$
when $\lambda\to 0$,
where the remainder term is locally uniform in $\hat x$ and with
$$
a_j(\hat x)={e^{\i\sigma_j(\hat x){\pi\over 4}}
\over
\sqrt{\vert K(k_j(\hat x))\vert } \vert \nabla \varepsilon(k_j(\hat x))\vert } 
$$
where $\sigma_j(\hat x)$ is the difference between the number of positive and negative principal curvatures of $\Sigma_E$ at $k_j(\hat x)$.
\end{lemma}
We first explain how to deduce Lemma~\ref{lem_Green} from previous Lemma~\ref{lem_Green_estimate}. Take a non singular unit vector $\hat x$ such that $\hat x_1>0$ and $\vert \hat x_1\vert >3 \vert \hat x_\ell\vert $ for $\ell=2, \ldots, d$ so that $L\hat x$ is in the set $C^{(L)}_{e_1}$.  By taking $\hat x$ in generic position and since $E\neq 2d$, we can always assume that $ \vert \nabla \varepsilon(k_j(\hat x))\vert$ are all distinct for $j=1, \ldots, m(\hat x)$
\footnote{ Indeed, if $k_j$ and $k_{j*}$ satisfies $| \nabla \varepsilon(k_j)| =| \nabla \varepsilon(k_{j*})| $ then  $\nabla \varepsilon(k_j)= \nabla \varepsilon(k_{j*})$  since they are both parallel and in the direction of $\hat x$. Hence, $\sin(k_j\cdot e_\ell)= \sin(k_{j*}\cdot e_\ell):= u_\ell, \; \forall \ell=1, \cdots, d$, and $\cos(k_j\cdot e_\ell)=\sigma_\ell \sqrt{1-u_\ell^2}$ and $\cos(k_{j*}\cdot e_\ell)=\sigma^*_\ell \sqrt{1-u_\ell^2}$ for $\sigma_\ell $ and $\sigma^*_\ell $ in $\{\pm 1\}$. Besides since $k_j, \; k_{j*}\in \Sigma_E$ then $\sum_{\ell=1}^d   \sigma_\ell \sqrt{1-u_\ell^2} = \sum_{\ell=1}^d   \sigma^*_\ell \sqrt{1-u_\ell^2} = 2d-E$. If we take $\hat x$ in generic position, hence $(u_\ell)$ generic, there is a  unique choice of signs $\sigma_\ell$, unless $2d-E=0$ in which case $(\sigma^*_\ell)=-(\sigma_\ell)$ and $k_{j*}=(\pi, \ldots, \pi)-k_j$ are dual solutions.
}. 
We can assume, without loss of generality, that $ \vert \nabla \varepsilon(k_1(\hat x))\vert > \vert \nabla \varepsilon(k_j(\hat x))\vert$   for $j=2, \ldots, m(\hat x)$. Taking $x=L\hat x$ with $L={K\over \lambda}$ for $K>1$, we see that we can take $K$ large enough so that
\begin{eqnarray*}
&&{\vert a_1(\hat x)\vert \over
\sqrt{\vert K(k_1(\hat x))\vert }} {e^{-\|x\|  \lambda   \over  \sqrt{2} \vert \nabla \varepsilon(k_1(\hat x))\vert}} 
= {\vert a_1(\hat x)\vert \over
\sqrt{\vert K(k_1(\hat x))\vert }} {e^{-K  \over \sqrt{2}  \vert \nabla \varepsilon(k_1(\hat x))\vert}} 
\\
&& > 2\l(\sum_{j=2}^{m(\hat x)} {\vert a_j(\hat x)\vert \over
\sqrt{\vert K(k_j(\hat x))\vert }} {e^{- K   \over  \sqrt{2} \vert \nabla \varepsilon(k_j(\hat x))\vert}}\r)
=2\l(\sum_{j=2}^{m(\hat x)} {\vert a_j(\hat x)\vert \over
\sqrt{\vert K(k_j(\hat x))\vert }} {e^{-\|x\|  \lambda   \over  \sqrt{2}   \vert \nabla \varepsilon(k_j(\hat x))\vert}}\r)
\end{eqnarray*}
and that it is true in a small neighborhood of $\hat x$ and $K$. In particular, it implies that there exists  $\epsilon>0$ such that, if $B(0,R)$ is the Euclidian ball with radius $R$, then for all $i\in B(0,L\epsilon)\cap \Z^d$ and $j\in L\hat x+ B(0,L\epsilon)\cap \Z^d$, we have, writing $y=L\hat x+j-i$, $K'={\| y \|\over \lambda}$,
\begin{eqnarray*}
\vert g_{{\i\over \sqrt{2}}\lambda -E}(i, j)\vert &= & \vert g_{{\i\over \sqrt{2}}\lambda -E}(0, y )\vert
 \ge {1\over 2} 
({\lambda\over K})^{d-1\over 2}{\vert a_1(\hat y)\vert\over 
\sqrt{\vert K(k_1(\hat y))\vert }} {e^{- K'  \over \sqrt{2}  \vert \nabla \varepsilon(k_1(\hat y))\vert}} 
\ge  C \lambda^{(d-1)/2}.
\end{eqnarray*}
since $(\hat y, K')$ are in a small neighbourhood around $(\hat x, K)$ and since $a_1(k_j(\hat y))$, $K(k_1(\hat y))$ and $\vert \nabla \varepsilon(k_1(\hat y))\vert$ are continuous. Hence, we have found for $\lambda$ small enough 
$$
B_0=B(0,L\epsilon)\cap \Z^d \subset C_{0}^{(L)} \;\hbox{ and }\; B_{e_1}=(B(0,L\epsilon)+L\hat x_1) \cap \Z^d \subset C_{e_1}^{(L)}
$$ 
which satisfy the requirement of the Lemma. 

By symmetry of the Green function, it means that we can find $\hat x_1, \ldots, \hat x_d$ for each base vector $e_1, \ldots e_d$ which satisfy the same conditions as above. Define then $\hat y_k:= \hat x_k-e_k$ for $k=1 \ldots, d$ and define for all $x\in \Z^d$:
$$
B_x:= B(0,L\epsilon) +Lx + L(\sum_{k=1}^d \hat y_k (x_k \bmod 2 )),
$$
which is included in $ C_{x}^{(L)}$ if we took $\hat y_k$ small enough. Then, $B_{x+e_k}$ is equal to the ball $B_x$ shifted by $L e_k\pm L \hat y_k$. Hence, they satisfy the conditions of Lemma~\ref{lem_Green}.
\end{proof}
 \begin{proof}[Proof of Lemma~\ref{lem_Green_estimate}]
 We come back to the Fourier representation \eqref{Fourier} for $g_{{\i\over \sqrt{2}}\lambda -E}(0,x)$.
We first localize around the Fermi surface.
We take $\delta>0$ small enough we denote $U_\delta=\{ k, \; |\varepsilon(k) - E|<\delta\}$, and set $\chi:\R\to [0,1]$, $C^\infty$, such that
$$
\chi(u)=\begin{cases}
\chi(u)=0,& \;\hbox{if $| u|\ge\delta$ }
\\
\chi({u})=1,&\;\hbox{if $| u|\le\delta/2$,}
\end{cases}
$$
 and we decompose the integral as
  \begin{eqnarray*}
&& \int_{\mathbb{T}^d } \frac{e^{i k \cdot {x}}}{\varepsilon(k) - E - \i\lambda } \, \mathrm{d}k 
 \\
 &=&\int_{\mathbb{T}^d} (1-\chi(\varepsilon(k)-E))\frac{e^{i k \cdot {x}}}{\varepsilon(k) - E - \i\lambda } \, \mathrm{d}k +   \int_{U_\delta } \chi(\varepsilon(k)-E) \frac{e^{i k \cdot {x}}}{\varepsilon(k) - E - \i\lambda } \, \mathrm{d}k 
 \end{eqnarray*}
By the non-stationary phase (or integration by part), since $| \varepsilon(k) - E|> \delta/2 >0$ when $(1-\chi(\varepsilon(k)-E))\neq 0$ the first term has rapid decay of order $\mathcal{O}\left(\|x\| ^{-N}\right)= \mathcal{O}\left(\lambda^{N}\right)$ for all $N$
uniformly in $\hat{x} \in \mathbb{S}^{d-1}$.

For $\delta$ small enough, $\nabla \varepsilon(k) \neq 0$ on $U_\delta$, and we apply the change of variables $k \mapsto (E', \sigma_{E'})$ where $\mathrm{d}\sigma_{E'}$ is the Euclidean surface measure on $\Sigma_{E'} = \{k \in \mathbb{T}^d \mid \varepsilon(k) = E'\}$ with $\mathrm{d}k = \frac{\mathrm{d}E' \, \mathrm{d}\sigma_{E'}(k)}{|\nabla \varepsilon(k)|}$, where $\mathrm{d}\sigma_{E'}$ is the Euclidean surface measure on $\Sigma_{E'} = \{k \in \mathbb{T}^d \mid \varepsilon(k) = E'\}$.

The second integral can be written as:
\begin{equation}
    g_{\i\lambda - E }(0,{x}) = \frac{1}{(2\pi)^d} \int_{E-\delta}^{E+\delta}
    \chi(E'-E) \frac{I(E', {x})}{E' - E - \i\lambda} \, \mathrm{d}E' + \mathcal{O}(\lambda^{N})
\end{equation}
where :
\begin{equation}
    I(E', {x}) = \int_{\Sigma_{E'}} e^{i k \cdot {x}} \frac{\mathrm{d}\sigma_{E'}(k)}{|\nabla \varepsilon(k)|}
\end{equation}

For $E' \in [E-\delta, E+\delta]$, the phase function $k \cdot {x}$ restricted to the surface $\Sigma_{E'}$ admits stationary points $k_j(E')$ when the unit normal ${n}(k) = \frac{\nabla \varepsilon(k)}{|\nabla \varepsilon(k)|}$ is aligned with $x$. Indeed, in this case, the phase is locally constant at first order on the Fermi manifold.
We denote by $(k^+_j(E'))_{i=1, \ldots m(\hat x)}$ (we forget here the dependence in $\hat x$ and focus on that of $E'$, so the notation is a bit different from the statement of the lemma) the stationary points which are aligned and in the direction of $n(k^+_j(\hat x))$ and $(k^-_j(E'))_{i=1, \ldots m(\hat x)}$ the ones which are aligned but in opposite direction to the normal (note that by symmetry of the Fermi curve we can take $k^-_j=-k^+_j$).

Since the vector $\hat x$ is non-singular, the Gaussian curvature does not vanish at $k_j^\pm(E')$ (at least for $\delta$ small enough), and the stationary phase theorem (see e.g. Theorem~9 page 30, \cite{Vainberg_book}) with parameters depending smoothly on $E'$ and $\hat{x}$ gives:
\begin{equation}
    I(E', {x}) = 
    \left( {2\pi \over \| x\|  }\right)^{\frac{d-1}{2}} \sum_{k^\pm_j} \frac{e^{i k^\pm_j(E')\cdot {x} + \i \frac{\pi}{4} \sigma^\pm_j}}{\vert \varepsilon(k^\pm_j(E'))\vert \sqrt{|K(k^\pm_j(E'))|}} + \mathcal{O}\left(\|x\|^{-\frac{d+1}{2}}\right)
\end{equation}
where $\sigma^\pm_j$ denotes the signature of the second fundamental form of $\Sigma_{E'}$ at $k^\pm_j$.

The last step is to integrate each term with respect to $E'$. When $\lambda\to 0$, we can do a classic contour integral. The term $ I(E', {x})$ is locally analytic in $E'$ and we restrict the integral to $\int_{E-\delta/2}^{E+\delta/2}$ to have $\chi(E'-E)=1$.
We now make a contour integral in the upper-half plane for $k^+_j$ or on the lower half-plane for $k^-_j$ to have a exponentially small remainder term. Hence, the stationary points $k^-_j$ give no residue, while the points $k^+_j$  give a residue at $E'=E+\i\lambda$.
Then, we make an expansion of $k^\pm_j(E') \cdot \hat{x}$ around $E' = E$. Differentiating the relation $\varepsilon(k^+_j(E')) = E'$ gives:
\begin{equation}
    \nabla \varepsilon(k^+_j) \cdot \frac{\mathrm{d}k^+_j}{\mathrm{d}E'} = 1 \implies {n}(k^+_j) \cdot \frac{\mathrm{d}k^+_j}{\mathrm{d}E'} = \frac{1}{\vert \varepsilon(k^+_j(E))\vert}
\end{equation}
Since ${n}(k^+_j) = + \hat{x}$, it follows that:
 $k^+_j(E') \cdot \hat{x} = k^+_j(E) \cdot \hat{x} + \frac{ 1}{\vert \varepsilon(k^+_j(E))\vert} (E' - E) + \mathcal{O}((E'-E)^2)$. This term is responsible for the exponential term $e^{-{\|x\| \lambda \over \sqrt{2} \vert\varepsilon(k_j(\hat x))\vert}}$, where we denotes $k_j(\hat x)$ for $k^+_j(E)$ in the statement. The part of the integral $\int_{\delta/ 2< |E'-E|<\delta}$ gives an extra term $\mathcal{O}\left(\| x\|^{-1}\right)$ by integration by part on $E'$.
 \end{proof}

  \begin{proof}[Proof of Theorem~\ref{thm_delocalization}~(ii)]
 The proof of the localization uses rather standard tools as the Combes-Thomas estimate to control the exponential decay of $W_{i,j}$ and the fractional moment method introduced in \cite{Aizenman_Molchanov,Aizenman94}. The fractional moment method was used to give a new proof of the localization of the $\H^{2|2}$ model in the nearest neighbor case in \cite{Collevecchio_Zeng21}. Our proof is rather close to that, however, we need to treat the case where the conductances are non local, which is done using Theorem~9.2 of \cite{Aizenman_Warzel}.
 
 Let us prove the exponential decay of the Green function.
 Let $\beta$ be the random potential with law $\nu_{\Z^d}^{rW}$ on $\Z^d$ and $\hat G_\beta$ be defined as in \eqref{Ghat}, which we simply write $\hat G$ for simplicity. As an almost sure limit of $\hat G^{(n)}$, $\hat G$ also admits a path expansion as in \eqref{G_paths}. Besides, $\hat G_{\beta}\le G_{\beta,\gamma}$ defined in \eqref{Ggamma},  which satisfies Proposition~\ref{law_gamma}.

 Take $y_0\in \Z^d$. Then $\E\l((\hat G_\beta(0,y)^s)\r)<\infty$ for $s<\demi$. By one-step expansion, we have for $x\neq y_0$,
 $$
 \hat G(x,y_0)=\hat G(x,x) \sum_{z\in \Z^d, z\neq x} r W_{x,z} \hat G^{\Z^d\setminus\{x\}}(z,y_0),
 $$
 where $\hat G^{\Z^d\setminus\{x\}}$ is the Green function of the operator restricted to $\Z^d\setminus\{x\}$, which can be defined by the path expansion \eqref{G_paths} for paths avoiding $x$. In particular, since $\hat G(z,y_0)$ is a sum on paths with positive terms we have 
 $$
 \hat G^{\Z^d\setminus\{x\}} (z,y_0)\le \hat G(z,y_0).
 $$
 Hence we write
 \beq
\E\l( \hat G(x,y_0)^s\r) &\le & \E\l(  G_{\beta,\gamma} (x,x)^s \sum_{z\in \Z^d} (rW)^s_{x,z} \hat G^{\Z^d\setminus\{x\}}(z,y_0)^s\r)
\\
&=& \E\l(  G_{\beta,\gamma} (x,x)^s\r)  \sum_{z\neq x} (rW)^s_{x,z} \E\l(  \hat G^{\Z^d\setminus\{x\}}(z,y_0)^s\r)
\\
&\le &
\E\l(   G_{\beta,\gamma} (x,x)^s\r)  \sum_{z\neq x} (rW)^s_{x,z} \E\l( \hat G (z,y_0)^s \r)
\eeq
since $ G_{\beta,\gamma}(x,x)$ is independent of $(\beta_y)_{y\neq x}$ and $\hat G^{\Z^d\setminus\{x\}}(z,y_0)^s$ measurable with respect to $(\beta_y)_{y\neq x}$.
Set now 
$$
a=\E\l(   G_{\beta,\gamma}(0,0)^s\r)=\E\l(   G_{\beta,\gamma}(x,x)^s\r) =\E\l(  \l({1\over 2\gamma}\r)^s\r) <\infty,
$$
and
$$
f(x)=\E\l( \hat G(x,y_0)^s\r).
$$
We have the following inequality
\begin{eqnarray} \label{ineq_f}
0\le f(x)\le g(x) + \sum_z K(x,z) f(z),\;\;\;\hbox{where}\;\;\; g(x)=a \delta_{x=y_0}
\end{eqnarray}
where 
$$
K(x,z)= \delta_{x\neq z} a (r W_{x,z})^s.
$$
We are in the context of Theorem~9.2 of \cite{Aizenman_Warzel}. Define $\tau(x)=\exp(+\mu\| x-y_0 \|)$. Since $W$ is exponentially decreasing by Combe-Thomas estimates (see e.g. \cite{Aizenman_Warzel} Theorem 10.5), if we assume that
\begin{equation}\label{sum_K}
\sum_{z} K(0,z) <1,
\end{equation}
then there is $\mu>0$ such that 
$$
b:=\sup_{x\in \Z^d}\sum_z { \tau(x)\over \tau(z)} K(x,z)\le \sup_{x\in \Z^d}\sum_z  \exp(\mu \| z-x\| ) K(x,z)=\sum_z  \exp(\mu \| z\| ) K(0,z) <1.
$$
Following the proof of Theorem~9.2 of \cite{Aizenman_Warzel}, it gives iterating \eqref{ineq_f} under condition \eqref{sum_K}
$$
f(x)\le \lim_{N\to \infty} \sum_{m=0}^N (K^m g)(x),
$$
and
$$
\sum_{x} \tau(x) f(x)\le \lim_{N\to \infty} \sum_{m=0}^N \sum_x \tau(x) (K^m g)(x)\le  
\lim_{N\to \infty} \sum_{m=0}^N b^m \sum_x \tau(x)  g(x)\le {a\over 1-b}.
$$
It implies that $\E\l ( (\hat G_\beta(x,y_0))^s\r )\le e^{-\mu \|x-y_0\|} {a\over 1-b}$, which gives exponential decay under condition \eqref{sum_K}. In particular, since $\hat G^{(n)}\le \hat G$ and by \eqref{Gbeta}, condition \eqref{sum_K} implies also exponential decay of $\E(\psi_n^s)$ and in particular that $\psi_\infty=0$ a.s.

We now have to check that the bound $\eqref{sum_K}$ is satisfied under the two assumptions of Theorem~\ref{thm_delocalization}~(ii). We apply Combes-Thomas estimate in the form of Theorem 10.5 \cite{Aizenman_Warzel} which we recall in our context: let $\delta:= \operatorname{dist}_\C(E+{\i\over \sqrt{2}}\lambda, \sigma(-\Delta)) $ and $S_\mu=2d(e^\mu -1)$ where $S_\mu $ is (10.31) of  \cite{Aizenman_Warzel} for $\mu>0$ and the nearest neighbor Laplacian. If $S_\mu<\delta$, then
$$
\vert g_{{\i\over \sqrt{2}}\lambda -E}(x,y)\vert \le {1\over \delta-S_\mu} \exp(-\mu \| x-y\|).
$$
Consider the first assumption of Theorem~\ref{thm_delocalization}~(ii)  and $E\not\in\sigma(-\Delta)$. Then $\delta> \operatorname{dist}_\C(E, \sigma(-\Delta)):= \delta'>0$. Take $\mu>0$ so that $S_\mu=\delta'/2$, it gives the bound
$$
\vert g_{{\i\over \sqrt{2}}\lambda -E}(x,y)\vert \le {2\over \delta'} \exp(-\mu \| x-y\|),
$$
which is independent of $\lambda$. Since 
$$
W_{x,y}=\lambda^2 \vert g_{{\i\over \sqrt{2}}\lambda -E}(x,y)\vert^2,
$$
taking $\lambda $ small enough we have that $rW$ satisfies the bound $\eqref{sum_K}$.

Consider now the second assumption. In this case we have $\lambda\le \sqrt{2} \delta$. Take $S_\mu=\delta/2$, we have
$$
W_{x,y}=\lambda^2 \vert g_{{\i\over \sqrt{2}}\lambda -E}(x,y)\vert^2\le 
{4 \lambda^2 \over \delta^2} \exp(-2\mu \| x-y\|)\le 8 \exp(-2\mu \| x-y\|).
$$
We can thus find $\Lambda(r)$ such that if $\delta\ge \Lambda(r)$, $\mu$ is large enough so that the bound \eqref{sum_K} is satisfied.
\end{proof}

\section{Discussion}\label{sec_discussion}
The main question we want to raise is wether the similarities between the two Brownian interpolations of the Anderson model with Cauchy potential and the $\H^{2|2}$ model is a mere coincidence or wether there is a hidden formal correspondence. Starting from the 
discussion of Section~\ref{sec_comparison},
we propose a conjecture.

Consider a finite graph as in Section~\ref{s_SDE_A}. We assume that this graph is vertex transitive, i.e. for all $x,y\in \cV$, there exists a permutations of the vertices sending $x$ to $y$ that leaves invariant the graph structure. Then, we want to consider the Anderson model with Cauchy potential with a uniform pinning $\epsilon>0$, i.e., for energy-disorder parameter $(E,\lambda)$ we are interested in the Schr\"odinger operator
\begin{eqnarray}\label{H_V_i_eps}
H_{V+{\i}{\epsilon\over \sqrt{2}}}= {1\over \lambda} (-\Delta -E)+V+{\i}{\epsilon\over \sqrt{2}}, \;\;\; G_{V+{\i}{\epsilon\over \sqrt{2}}}=\l ( H_{V+{\i}{\epsilon\over \sqrt{2}}} \r )^{-1}
\end{eqnarray}
where $(V_v)_{x\in \cV}$ are i.i.d. with Cauchy law with scale ${1\over \sqrt{2}}$. Then, for $0<p<1$, we consider
\begin{eqnarray}\label{IA}
I_A=
\E \l( \l | \sqrt{2} (\Im(G_{V+{\i}{\epsilon\over \sqrt{2}}}(i_0,i_0)))\r | ^p\r),
\end{eqnarray}
which does not depend on $i_0\in \cV$ by vertex transitivity.

On the $\H^{2|2}$ side, for $r>0$, we consider the finite set $\cV$ with conductances $rW$ where $W_{x,y}$ is given, for $x,y\in \cV$, by
\begin{eqnarray}\label{W_again_i_eps}
W_{x,y}=\l \vert \l({1\over \lambda} (-\Delta -E)+{\i\over \sqrt{2}} +\i{\epsilon\over \sqrt{2}}  \r)^{-1}_{x,y}\r\vert^2,
\end{eqnarray}
and a boundary term $\eta$, given for all $x\in \cV$ by
\begin{eqnarray}\label{eta_i_eps}
\eta_x
=\epsilon \sum_{y \in \cV} W_{x,y},
\end{eqnarray}
which is constant on $\cV$ by vertex transitivity. Consider now the random potential $(\beta_x)_{x\in \cV}$ with distribution $\nu_{\cV}^{rW,r\eta}$ defined in Theorem~\ref{prop:nuwtheta}, Section~\ref{sec_def_beta}, and we set 
$$
\psi_x=(G_{\beta} r \eta)_x, \;\;\; G_{\beta}=(\beta-rW)^{-1}.
$$
Note that $\psi_x$ is also the first marginal of the $\H^{2|2}$ model with uniform pinning as in \cite{DSZ10} (where there $W$ is nearest neighbor conductances on a box of $\Z^d$).
Then, we set
\begin{eqnarray}\label{IS}
I_{r,S} = \E_{\nu^{rW,r\eta}} \l((\eta_{i_0} +(W  \psi)_{i_0} )^p\r) .
\end{eqnarray}
which again does not depend on the choice of $i_0\in \cV$ by vertex transitivity.

In Section~\ref{sec_relation_SDE} we explain how this construction follows naturally from the comparison between the stochastic differential equations for the $\H^{2|2}$ and the Anderson model with frozen vertices investigated in sections~\ref{sec_SDE_A_bord} and  \ref{sec_SDE_S_bord}. The quantities $I_A$ and $I_{r,S}$ are the fractional moments of the limit values of the natural martingales that appear in both cases.

We think that the two models are related via the moments of the limit values of these martingales. We propose a conjecture.
\begin{conjecture}\label{conjecture}
Let $0<p\le {1\over 2}$. There exists $0< r_0(p)<1$, such that for any finite vertex transitive connected graph $\mathcal{G}$ with uniform pinning as above, and any $(E,\lambda)\in \R\times \R^*_+$, 
\begin{eqnarray*}
I_{r_0(p),S}\le I_A\le I_{1,S}.
\end{eqnarray*}
\end{conjecture} 

The conjecture comes from the similarities we identified between the Brownian interpolations of the two models and several numerical simulations we performed on $\Z^d$, $d=1,2,3$, which showed that the two models were rather similar for all values of $(E,\lambda)$ for the same choice of $r$. However, our first version of the conjecture for all graphs with a frozen subset was shown to be wrong when we checked it with ChatGPT6.
The version above seems to be robust to numerical AI-assisted experiments.

Besides, we propose in next Section~\ref{sec_one_point} a stronger version of the conjecture that gives an explicite value for the lower bound. In fact, it seems that the case of one point with a pinning $\epsilon$ is extremal for the lower bound for $0<p\le  \demi$. 

Note that proving the result for any $p$ is sufficient for applications. The choice $p\le \demi$ comes from the fact that at very small pinning the Green function $G_\beta$ for $\H^{2|2}$ and the corresponding $| G_V|^2$ for Anderson have exploding $p$-moments for $p\ge 1/2$. We will later update the paper with a few special cases where we can compare the phases or moments of the two models. 

\subsection{Implications it would have on the delocalization for Anderson with Cauchy potential on $\Z^d$, $d\ge 3$.} 
Let us explain the implications this conjecture would have in the case of $\Z^d$. We do not detail all arguments since the conjecture is not proved.

Consider the torus of size $2N$: $\T_N=(\Z / 2N \Z)^d\eqsim [-N,N[^d \cap \Z^d$.  Fix the energy and disorder parameters $(E,\lambda)$.  Consider now the Anderson model with Cauchy potential and uniform pinning $\epsilon>0$. The Schr\"odinger operator is then
$$
H_{N,V+\i{\epsilon\over \sqrt{2}}}={1\over \lambda}(-\Delta_{\T_N}-E) + V+ \i{\epsilon\over \sqrt{2}}, \;\;\, G_{N,V+\i{\epsilon\over \sqrt{2}}}=(H_{N,V+\i{\epsilon\over \sqrt{2}}})^{-1},
$$
where $\Delta_{\T_N}$ is the Laplace operator on $\T_N$. Hence, we have, with definition \eqref{IA},
$$
I^{(N)}_A=\E \l(  \l | \sqrt{2} \Im\l( G_{N,V+\i{\epsilon\over \sqrt{2}}}(0,0)\r) \r |^p\r).
$$

On the $\H^{2|2}$ side the conductances are given by $r W^{(N,\epsilon)}$ where, for $x,y\in \T_N$,
$ W^{(N,\epsilon)}_{x,y}$ is given by \eqref{W_again_i_eps}, with the discrete Laplacian $\Delta_{\T_N}$ on the torus,
and boundary condition $r \eta^{(N,\epsilon)}$, where
$
\eta^{(N,\epsilon)}_{x}= \epsilon \sum_{y\in \T_N} W_{x,y}^{(N,\epsilon)}.
$
Then, with definition \eqref{IS}
$$I^{(N)}_{r,S}=
\E
 \l ( \l( \eta_0^{(N,\epsilon)}+ (W^{(N,\epsilon)} \psi^{(N,\epsilon))})_0\r)^p\r),
$$
where 
$ \psi^{(N,\epsilon))} = G_{\beta^{(N,\epsilon)}} r\eta^{(N,\epsilon)}$, and where $\beta^{(N,\epsilon)}$ follows the law $\nu_{\T_n}^{rW^{(N,\epsilon)}, r\eta^{(N,\epsilon)}}$ and $G_{\beta^{(N,\epsilon)}}=(\beta^{(N,\epsilon)} - r W^{(N,\epsilon)})^{-1}$.

Then, letting $N$ go to $\infty$ 
we would get
$$
\lim_{N\to \infty} I^{(N)}_{r,S}= \E_{\nu_{\Z^d}^{W^{(\epsilon)},\eta^{(\epsilon)}}}
\l ( \l( W^{(\epsilon)} \psi^{(\infty))}\r)_0^p\r)\ge \E_{\nu_{\Z^d}^{W^{(\epsilon)},0}}
\l ( \l( W^{(\epsilon)} \psi^{(\infty))}\r)_0^p\r) ,
$$
where $W^{(\epsilon)}_{x,y}$ and $\eta^{(\epsilon)}$ are defined in \eqref{W_again_i_eps} and \eqref{eta_i_eps} with uniform pinning $\epsilon>0$ on $\Z^d$. The second inequality is a consequence of monotonicity.
An attentive reading of the proof of Theorem~\ref{thm_delocalization}~(i) shows that the lower bound on $\E((\psi^{(\infty)})^p)$ can be proved for conductances $W^{(\epsilon)}$ instead of $W$, uniformly in $\epsilon$ for $\epsilon$ small enough. It implies that for all $E\in (0,4d)\setminus \sss_0$, for $\lambda<\lambda(E,r_0(p))$
$$
\liminf_{\epsilon\to 0}\lim_{N\to \infty} I^{(N)}_{r_0(p),S}>0
$$
locally unformly in $E$.
%

On the Anderson side, we have
$$
\lim_{N\to \infty} I_A= \E \l(  \l | \sqrt{2} \Im\l( G_{\Z^d,V+\i\epsilon{1\over \sqrt{2}}}(0,0)\r) \r |^p\r).
$$
where $G_{\Z^d,V+\i\epsilon{1\over \sqrt{2}}}$ is the Schrodinger operator on $\Z^d$ with uniform pinning $\epsilon$. If the conjecture were true then it would imply
$$
\liminf_{\epsilon\to 0} \E \l(  \l | \sqrt{2} \Im\l( G_{\Z^d,V+\i{\epsilon\over \sqrt{2}}}(0,0)\r) \r |^p\r)
>0.
$$
Hence it would imply existence of an absolutely continuous part in the spectrum of the Anderson model with Cauchy potential in a small neighborhood around spectral value $E$, using \cite{Jaksik} Theorem~3.16 together with an elementary exchange of limit and integral.

Note however, that even if conjectures~\ref{conjecture} and \ref{conjecture_2points} below are true it is not sure that they are easier to prove than delocalization.

\subsection{How conjecture~\ref{conjecture} is related to the stochastic differential equations of Sections \ref{s_SDE_A} and \ref{s_SDE_S}}\label{sec_relation_SDE}

In order to fit Conjecture~\ref{conjecture} in the context of Section~\ref{sec_SDE_A_bord}, we consider now a copy $\tilde \cV$ of the set $\cV$ and define $\tilde \ggg=(\cV\cup \tilde \cV,\tilde \cE)$ where we keep the edges of the graph $\ggg$ between the vertices $\cV$ and add exactly one edge between $x$ and its copy $\tilde x$ in $\tilde \cV$. We put the conductance $\Delta_{x,\tilde x}=\lambda \sqrt{{\epsilon\over{2}}}$ for all $x\in \cV$ and a diagonal shift in the Laplacian $\Delta_{\tilde x,\tilde x}=-E$. Then we consider the vertices of $\tilde \cV$ as frozen, i.e. we put a potential $V_{\tilde x}={\i\over \sqrt{2}}$ on the vertices $\tilde x\in \tilde \cV$.  It means that the Schur complement on $\cV$ of the full Schr\"odinger operator on $\cV\cup \tilde \cV$ has exactly the form~\eqref{H_V_i_eps}. Besides, by Schur complement, if we compute the Green function on the full graph $\tilde\ggg$ with frozen potential $\i/\sqrt{2}$ on all vertices of $\cV\cup\tilde\cV$ and set
 $$ W:= 
 | G_{\cV\cup\tilde\cV, {\i\over \sqrt{2}}} |^2,
 $$
 then we get that $W_{x,y}$ is given by \eqref{W_again_i_eps} for $x,y\in \cV$ active vertices and $W_{x,\tilde y}=\epsilon W_{x,y}$ for $x\in \cV$ and $\tilde y\in \tilde \cV$. Since the vertices of $\tilde\cV$ are frozen, by Section~\ref{sec_SDE_S_bord}, the natural $\H^{2|2}$ model to consider  is the model with conductances $rW$ on $\cV\times \cV$ and boundary term $(r\eta_x)_{x\in \cV}$ given by \eqref{eta_i_eps}.
 
 Now we consider the stochastic processes $(W^A(t),X^A(t))$ and $(W^S(t),X^S(t))$ defined in Section~\ref{sec_SDE_A_bord} and \ref{sec_SDE_S_bord}, where $W^{A/S}(t)$ is defined on $(\cV\cup \tilde \cV)\times( \cV\cup \tilde \cV)$ and $X^{A/S}(t)$ on $\cV\cup \tilde \cV$.  The vertices of $\tilde \cV$ are frozen and $X^{A/S}_{\tilde \cV}(t)=1$ $\forall t>0$. 
 
 For all $x\in \cV$, $(W^{S/A}(t) X^{S/A}(t))_{x} $ are martingales. By Ward identity \eqref{Ward}, we have for $x\in \cV$,
$$
(W^{A}(t) X^{A}(t))_{x}   = - \sqrt{2} (\Im(G_{V(t)+{\i}{\epsilon\over \sqrt{2}}}(x,x))),
$$
where $G_{V(t)+{\i}{\epsilon\over \sqrt{2}}}$ is the Green function on $\cV$ defined in \eqref{H_V_i_eps}.
Hence, $I_A$ is the fractional moment of the limit value of the martingale: $I_A=\E\l( (W^{A}(\infty) X^{A}(\infty))^p_{i_0}\r) $ for any $i_0\in \cV$.

On the $\H^{2|2}$ side, the process $(W^S(t),X^S(t))$ is a diffusion. Starting from initial values $W^S(0)=rW^A(0)= rW$, we get for $x\in \cV$,
$$
(W^S(\infty)X^S(\infty))_x=(W^S(0) (K_{T(\infty)})^{-1} X^S(\infty))_x=  r (\eta_x+(W \psi)_x).
$$
Hence, $I_{r,S}={1\over r^p} \E\l ( (W^S(\infty)X^S(\infty))^p_x \r)$ is the fractional moment of the limit value of the martingale, with natural rescaling by ${1\over r}$ since $ (W^S(0)X^S(0))_x= r (W^A(0)X^A(0))_x$.

\subsection{Heuristic}\label{sec_Heuristic}
We try to give an heuristic to substantiate this conjecture, remaining very informal. Even if it does not reflect the complexity of the problem (in particular, it would apply to cases were the conjecture is false),
it highligths some cancellations in the generators of the two processes. We consider the processes $(W^{A/S}(t), X^{A/S}(t))$ defined in section~\ref{sec_SDE_A_bord} and \ref{sec_SDE_S_bord}, with active subset of vertices $\cV$ and frozen subset of vertices $\tilde \cV$. (Mind that the notations are slightly different from   section~\ref{sec_SDE_A_bord} and \ref{sec_SDE_S_bord} where the full set was $\cV$ and the active set $\cV'\subset \cV$).
We consider the generators $L^A$ of the process $(W^A(t),X^A(t))$, and $L^S$ of the process $(W^S(t),X^S(t))$, and $P^A_t$ and $P^S_t$ the associated semi-groups. Note that since $(W^A(t),X^A(t))$ is a priori not Markov (it is rather $(X(t),Y(t))$ which is Markov), the differential operator $L^A$ will contain terms which are not measurable with respect to $(W,X)$ and $P^A_t(\phi)$ means the expectation of $\phi(W^A(t),X^A(t))$ and $L^A$ satisfies $\partial_t P_t^A(\phi)=P_t^A L^A \phi$. In  particular, as we will compose operators, the only licit composition is $P^A_t\circ P^S_s$ and not the reverse.

For $(W,X)$ we set $\phi(W,X)=(WX)_{i_0}^p$ for a fixed vertex $i_0$. With the notations of Section~\ref{sec_relation_SDE} above, we have $I_A= (P^A_\infty \phi)$ and   $I_{r,S}={1\over r^p} (P^S_\infty \phi)(r W,X=1)$.   
To simplify notations we will write
$$
F(W,X)=(P^S_\infty \phi)(W,X)
$$
and for $r>0$ we set
$$
F_r(W,X)={1\over r^p} F(rW,X),
$$
 and also $\tau_r$ the application $f\circ\tau_r(W,X):= f(r W,X)$.  Note that for $t=\infty$ we have $X_i(\infty)=\delta_{i\in \tilde \cV}$. In particular we have $F_r(W^A(\infty),X^A(\infty))={1\over r^p}\phi(rW^A(\infty),\delta_{i\in \tilde\cV})=\phi(W^A(\infty),X^A(\infty))$. Hence we have
 $$
 P^A_\infty(F_r)=P^A_\infty(\phi).
 $$
 
 The idea would be to use the structure of $L^A$ and $L^S$, where there are several coincidences due to similarity of their stochastic equations.

We write 
\begin{eqnarray}
\label{int_Pt}
I_A-I_{r,S}=(P^A_\infty \phi) - F_{r} (W,X=1)
&=& P^A_\infty \l (F_{r}\r) - F_{r}(W^A(0),X^A(0))
\\
\nonumber
&=&
\int_0^\infty {\partial\over \partial t} P^A_t\circ \l(F_{r} \r) dt
\\&=&
\nonumber
\int_0^\infty P^A_t\l(L^A(F_{r})\r)  dt
\end{eqnarray}
We are thus interested in the sign of \eqref{int_Pt}. We look below at the sign of $L^A(F_r)$ in some approximations for small and large $r$.

From \eqref{SDE_W_S}, direct computation of $L^S$ gives
$$
L^S = B+ \demi\sum_{k\in \cV} X_k^2 {\partial^2_{X_k^2}},
$$
with
$$
B= \sum_{i,j} \sum_{k\in \cV} (W_{i,k} X_k^2 W_{k,j}) \partial_{W_{i,j}} - \sum_{k\in \cV} X_k^2 (W X)_k \partial_{X_k}.
$$
Note that by Feynman-Kac, we have $L^S(F)=0$. Hence, it implies that $B(F)=- \demi\sum_{k\in \cV} X_k^2 {\partial^2_{X_k^2}} F$. If now we consider $F_r$, then we have $((B(F))\circ\tau_r)(W,X)=r^{p+1} B(F_r)(W,X)$. It implies in particular that 
\begin{eqnarray}\label{egal_Fr}
L^S_r(F_r)=0, \;\;\;\hbox{where } L^S_r=r B+ \demi\sum_{k\in \cV} X_k^2 {\partial^2_{X_k^2}}
\end{eqnarray}
From \eqref{eq_W_A}, we get
$$
L^A(F_r) = B(F_r)+ \tilde L^A(F_r),
$$
with
\begin{eqnarray*}
\tilde L^A&=& \sum_{i,j}\sum_{i',j'} \sum_{k\in \cV} X_k^2 \overline\sigma^{i',j'}_k \sigma^{i,j}_k \partial^2_{W_{i,j}W_{i',j'}} -  \sum_{i,j}\sum_{k\in \cV} \sqrt{2} \hbox{Im} (\sigma^{i,j}_k) X_k^2   \partial^2_{W_{i,j},X_k}
\\
&& + \demi\sum_{k\in \cV} X_k^2 {\partial^2_{X_k^2}} + \sum_{k\in \cV} X_k^2 (W X)_k \partial_{X_k}
\eeq
where $ \sigma^{i,j}_k$ is defined in \eqref{sijk}. The similarities are not obvious on these formula, they are more clear when we look at the Gauge variables $\tilde W_{x,y}$ since there are extra cancellation between the terms $\hbox{Im} (\sigma^{i,j}_k)X_k^2   \partial^2_{W_{i,j},X_k} $ when $k=i$ and the term $(X^2 W X)_k \partial_{X_k}$. However, we can look at special functions on which coincidences appear.

Take first the function $\phi(W,X)=(WX)_{i_0}^p$. In the very delocalized case, since $W(t)X(t)$ is a martingale and since there are only weak fluctuations of the function $\psi$, we have $F(W,X)\sim \phi(W,X)$. In particular, for $r$ very large, there is weak disorder on the $\H^{2|2}$ side, and $F_r(W,X)\sim \phi(W,X)$.

Let us first compute $\tilde L^A(F_r)$ in this limiting case. If we set
$$
\Sigma_k=\left(i\sum_{i,j} \sigma^{i,j}_k {\partial\over\partial W_{i,j}}\right) +{1\over \sqrt{2}} {\partial\over\partial X_k},
$$
we have,
$$
\tilde L^A((WX)_{i_0}^p)=
p(p-1) (WX)_{i_0}^{p-2} \sum_{k\in \cV} X_k^2 \overline\Sigma_k ((WX)_{i_0})  \Sigma_k ((WX)_{i_0}).
$$
Indeed, the first order derivative $\sum_{k\in \cV} (X^2WX)_k {\partial\over \partial X_k}$ acting on a single term $(WX)_{i_0}$ cancels with the action of  $\sum_{i,j}\sum_{k\in \cV} \sqrt{2} \hbox{Im} (\sigma^{i,j}_k) X_k^2   \partial_{W_{i,j},X_k}$ since only the terms $i=i_0$ and $k=j$ contribute, and thanks to the Ward identity \eqref{Ward}, $\sqrt{2} \sum_j \hbox{Im} (\sigma^{k,j}_k)= W_{k,j} (WX)_{k}$.

 As $i{ \sqrt 2} X=H- \overline H$, it gives
 \begin{eqnarray*}
i \sum_j \sigma^{i_0,j}_k X_j &=& - i \sum_j \overline{G(i_0,j)}  G(i_0,k)G(k,j) X_j= - G(i_0,k) {1\over \sqrt 2}  \left( G(H-\overline H)\overline G\right)_{i_0,k}
\\
&=& i { \sqrt 2}  G(i_0,k)  \hbox{Im}(G(i_0,k))
\end{eqnarray*} 
Hence,
\begin{eqnarray*}
\Sigma_k ((WX)_{i_0})&= & { \sqrt 2}  G(i_0,k)  \hbox{Im}(G(i_0,k)) +{1\over \sqrt 2} W_{i_0,k}
\\
&= &  G(i_0,k)\left( i ({ \sqrt 2}-{1\over \sqrt 2})    \hbox{Im}(G(i_0,k)) +{1\over \sqrt 2}  \hbox{Re}(G(i_0,k)) \right)
\\
&= &  {1\over \sqrt 2} G(i_0,k)^2.
\end{eqnarray*}
Finally, we have,
$$
\overline\Sigma_k((WX)_{i_0})  \Sigma_k ((WX)_{i_0})={1\over  2} \vert G(i_0,k)^2\vert ^2= {1\over  2} W_{i_0,k}^2,
$$
and
$$
\tilde L^A((WX)_{i_0}^p)=p(p-1) \left( \sum_k (W_{i_0,k} X_k)^2\right)(WX)_{i_0}^{p-2}=\demi\sum_k X_k^2 {\partial^2_{X_k^2}}((WX)_{i_0}^p).
$$
Now, we look at the term $B(F_r)$. Formally, we have $B((WX)_{i_0}^p)=0$, but it is incompatible with \eqref{egal_Fr}, and in order to take into account small terms in ${1\over r}$, we apply \eqref{egal_Fr} and write $B(F_r)=-{1\over r} (\demi\sum_k X_k^2 {\partial^2_{X_k^2}}(F_r))$, and we get for $r$ large,
$$
L^A(F_r))\sim (1-{1\over r}) p(p-1) \left( \sum_k (W_{i_0,k} X_k)^2\right)(WX)_{i_0}^{p-2}.
$$
which is negative for $r>1$  and $p<1$, hence, we see a upper barrier at $r=1$ in this approximation.

Consider now the limiting case $r\to 0$, in this case we have $F_r\sim (\sum_{j\in \tilde\cV} W_{i_0,j})^p$ since it is the very localized limit and almost no interaction is created during the interpolation and $X_x(\infty)=\delta_{x\in \tilde\cV}$. In this case $\partial^2_{X_k^2}(\sum_{j\in \tilde\cV} W_{i_0,j})^p=0$ and using \eqref{egal_Fr} we replace the degenerate term $\demi \sum X^2_k \partial^2_{X_k^2} (F_r)$ by $-rB(F_r)$ to take into account small terms in $r$. It gives, 
\begin{eqnarray*}
&&L^A( F_r )
\\
&\sim & (1-r)B((\sum_{j\in \tilde\cV} W_{i_0,j})^p)+ \sum_{i,j}\sum_{i',j'} \sum_{k\in \cV} X_k^2 \overline\sigma^{i',j'}_k \sigma^{i,j}_k \partial_{W_{i,j}W_{i',j'}}  ((\sum_{j\in \tilde\cV} W_{i_0,j})^p)
\\
&=&
p(1-r)(\sum_{j\in \tilde\cV} W_{i_0,j})^{p-1} \sum_{k}\sum_{j\in \tilde\cV} W_{i_0,k} X_k^2 W_{k,j}+ p(p-1) (\sum_{j\in \tilde\cV} W_{i_0,j})^{p-2} \sum_{k\in \cV} X_k^2 \l | \sum_{j\in \tilde\cV}  \sigma^{i_0,j}_k\r |^2 
\end{eqnarray*}
We have now by Cauchy-Schwarz inequality  and the expression of $\sigma_k^{i,j}$, \eqref{sijk},
$$
\l | \sum_{j\in \tilde\cV}  \sigma^{i_0,j}_k\r |^2 \le  \l ( \sum_{j\in \tilde \cV} \l | G(i_0,j)\r |^2 \r )  \l ( \sum_{j\in \tilde \cV} \l|G(i_0,k) G(k,j)\r |^2 \r )=  \l ( \sum_{j\in \tilde \cV} W_{i_0,j} \r )   \l ( \sum_{j\in \tilde \cV} W_{i_0,k} W_{k,j} \r )
$$
Hence, the last expression above gives since $p<1$,
\beq
L^A( F_r )&\gtrsim &
p(p-r) (\sum_{j\in \tilde \cV} W_{i_0,j})^{p-1} \sum_{k} \sum_{j\in \tilde \cV}   W_{i_0,k} X_k^2 W_{k,j}
\end{eqnarray*}
which is positive for $r<p$. Hence, we see a lower barrier at $r=p$ in this approximation.

Note that the heuristic above is by far not rigorous and even if the heuristic barriers at $p$ and 1 are uniform in the graph and parameters, the approximations $F_r\sim (WX)_{i_0}^p$ for $r$ large or $F_r\sim (\sum_{j\in \tilde \cV} W_{i_0,j})^p$ for $r$ small are absolutely not uniform. Note that the largest lower bound is not $p$ but in general smaller (see in particular the one-point graph in Section~\ref{sec_one_point}).

\section{The one point case with pinning and a more precise conjecture}\label{sec_one_point}

Consider the graph $o-\tilde o$, where $o$ is active and $\tilde o$ is frozen, with pinning $\epsilon>0$. We use the definitions of $I_A$, $I_{r,S}$, $W$ and $\eta$ in \eqref{IA}--\eqref{IS}. For $0<p<1$, set
\begin{eqnarray}\label{one_point_constants}
a_p&=&p-{\Gamma(p+\demi)\over\sqrt{\pi}\Gamma(p)},\nonumber\\
b_p(r)&=&{\sqrt r\over 2\sqrt\pi}\int_0^\infty
e^{-ru}u^{-3/2}\bigl((1+u)^p-1\bigr)\,du,\quad r>0.
\end{eqnarray}
The function $b_p$ is strictly increasing from $0$ to $p/2$, and $0<a_p<p/2$. In particular, the equation $b_p(r)=a_p$ has a unique solution, which we denote by $r_*(p)$. We have the following numerical estimates
$$
r_*(1/4)=0.03558366966675624\ldots,\qquad
r_*(1/2)=0.12052713163903786\ldots\,.
$$

\begin{proposition}\label{prop_2_points}
For the one point graph and all $0<p<1$, $\epsilon>0$, $E\in\R$ and $\lambda>0$,
\begin{eqnarray}\label{encadrement_o_delta}
I_{r_*(p),S}< I_A < I_{p,S}.
\end{eqnarray}
Moreover, $r_*(p)$ is the largest $r$ for which the lower inequality $I_{r,S}\le I_A$ holds for all $\epsilon$, $E$ and $\lambda$. 
\end{proposition}

Besides, numerical simulations 
seem to indicate that the two points graph gives the worst lower bound, which suggests the following more precise conjecture for the lower bound.
\begin{conjecture}\label{conjecture_2points}
For $0<p\le \demi$, the lower bound in Conjecture~\ref{conjecture} holds with $r_0(p)=r_*(p)$.
\end{conjecture}
Note that the upper constant $p$ in \eqref{encadrement_o_delta} concerns only the one point graph; the upper bound in Conjecture~\ref{conjecture} remains stated with $r=1$.
\begin{proof}[Proof of Proposition~\ref{prop_2_points}]
We could give a proof of a uniform control as in \eqref{encadrement_o_delta} that follows the Heuristic of Section~\ref{sec_Heuristic}. However, it does not give the optimal lower bound. In fact the largest $r$ that satisfies the punctual inequality $P_t^A (L^A F_r) >0$ for all $t>0$, is not the optimal one for the lower bound of Proposition~\ref{prop_2_points}. We propose a proof found by ChatGPT6 based on the explicit expression of $I_A$ and $I_{S,r}$ and which gives the optimal lower bound.  

Consider the graph with vertices $\{o,\tilde o\}$ where $o$ is active and $\tilde o$ frozen. Then on the full set $\{o,\tilde o\}$ we have, see Section~\ref{sec_relation_SDE},
$$
H_{V}=
{1\over \sqrt{2}} \left(
\begin{matrix}
 \sqrt{2} V_o-{\sqrt{2} E\over \lambda} & -{\sqrt{{\epsilon}}}
\\
-{\sqrt{{\epsilon}}}  & {\i}
\end{matrix} \right), 
\;\;\; 
H_{\i\over \sqrt{2}}=
{1\over \sqrt{2}} \left(
\begin{matrix}
 {\i}  -{\sqrt{2} E\over \lambda} & -{\sqrt{{\epsilon}}}
\\
-{\sqrt{{\epsilon}}}  & {\i} 
\end{matrix} \right).
 $$
Set $a=-\sqrt{2}E/\lambda$. The one point versions of \eqref{W_again_i_eps} and \eqref{eta_i_eps} give
$$
W_{o,o}=\l | (H_{\i\over \sqrt{2}}^{-1})_{o,o}\r |^2={2\over a^2+(1+\epsilon)^2 },\qquad \eta_o=\epsilon W_{o,o}={2\epsilon\over  a^2+(1+\epsilon)^2}.
$$
If $V=\sqrt{2}V_o$ is a standard Cauchy random variable, then
\begin{eqnarray}
I_A&=& \E\l( \l( - \sqrt{2} \Im ( ( H_{V}^{-1})_{o,o}) \r)^p \r) = 2^p\E\l(\l({\epsilon\over (V+a)^2+\epsilon^2}\r)^p\r),\nonumber\\
I_{r,S}&=&W_{o,o}^p\E\l((\epsilon+\psi_o)^p\r).\label{one_point_moments}
\end{eqnarray}
For the one point graph, since $\beta_o\sim\nu_{\{o\}}^{rW_{oo},r\eta_o}$, $\psi_o={1\over \beta_o-rW_{o,o}} r\eta_o$ 
has Laplace transform given by
$$
\E(e^{-s\psi_o})=
\exp\l[-r\eta_o\l(\sqrt{1+{2s\over r\eta_o}}-1\r)\r],\qquad s\ge0.
$$

We first express $I_A$ and $I_{S,r}$ through the same elementary function. Write
$$
x={1\over 1+\epsilon},\qquad R={a^2+(1+\epsilon)^2 \over(1+\epsilon)^2}\ge1,
\qquad \chi(u)={2\over1+\sqrt{1+u}}.
$$
Let $U$ and $T$ be positive random variables with respective densities on $(0,\infty)$
$$
{\sin(\pi p)\over\pi}{u^{-p}\over1+u}\,du,
\qquad {t^{-p}e^{-t}\over\Gamma(1-p)}\,dt.
$$
The normalization of the first density comes from the beta-prime integral and the duality of Gamma function: $\int_{0}^\infty {s^{-p}\over {1+s}}ds=\Gamma(p)\Gamma(1-p)={\pi\over \sin(\pi p)}$. 
Thus $U$ has the law beta-prime with parameters $(1-p,p)$, and $T$ has a gamma distribution with parameter $(1-p,1)$. 
The lemma below gives a common representation for $I_A$ and $I_{s,r}$.
\begin{lemma}\label{lem_common_rep}
\begin{eqnarray}
{I_A\over\eta_o^p}&=&\E\l((1-x\chi(RU))^{-p}\r),\nonumber\\
{I_{r,S}\over\eta_o^p}&=&\E\l((1-x\chi(RT/r))^{-p}\r).\label{one_point_common_rep}
\end{eqnarray}
\end{lemma}
\begin{proof}
For the first identity, we use the beta-prime integral to write
$$
q^{-p}={\sin(\pi p)\over\pi}\int_0^\infty {s^{-p}\over q+s}\,ds.
$$
Hence,
\begin{eqnarray*}
{I_A\over(2\epsilon)^p}&=&
 { \sin(\pi p)\over\pi} \int_0^\infty\int_0^\infty { s^{-p}\over (v+a)^2+\epsilon^2+ s} {1\over 1+ v^2} ds dv 
 \\
 &=& {\sin(\pi p)\over\pi} \int_0^\infty  {s^{-p}\over \sqrt{\epsilon^2+s} } \l( \int_0^\infty { \sqrt{\epsilon^2+s} \over (-a-v)^2+ ({\epsilon^2+s}) } {1\over 1+ v^2}  dv\r) ds 
 \end{eqnarray*}
 We recognize the convolution of a Cauchy with parameter 1 and a Cauchy with parameter $\sqrt{\epsilon^2+s}$, which is a Cauchy of parameter $(1+\sqrt{\epsilon^2+s})$. Hence it gives, changing to variable $b=\sqrt{\epsilon^2+s}$ in the second equality,
 \begin{eqnarray*}
{I_A\over(2\epsilon)^p}
&=& { \sin(\pi p)\over\pi} \int_0^\infty  {s^{-p}\over \sqrt{\epsilon^2+s} }  {1+\sqrt{\epsilon^2+s} \over a^2+(1+\sqrt{\epsilon^2+s} )^2} ds 
\\
&=&
{2\sin(\pi p)\over\pi}\int_\epsilon^\infty
(b^2-\epsilon^2)^{-p}{1+b\over a^2+(1+b)^2}\,db.
\end{eqnarray*}
Now substitute $u={(1+b)^2-(1+\epsilon)^2 \over a^2+(1+\epsilon)^2}$. A direct calculation gives
$$
{a^2+(1+b)^2 \over a^2+(1+\epsilon)^2}=(1+u),\qquad
{b^2-\epsilon^2\over a^2+(1+\epsilon)^2}=u[1-x\chi(Ru)],
$$
and proves the first formula in \eqref{one_point_common_rep} using $du=2(1+b)/(a^2+(1+\epsilon)^2)db$ and $\eta_o=\epsilon W_{oo}=2\epsilon/(a^2+(1+\epsilon)^2)$.

For the second identity, write
$$
\E(e^{-s(1+\psi_o/\epsilon)})=e^{-H(s)},\qquad
H(s)=s+{2r\epsilon\over D}
\l(\sqrt{1+{Ds\over r\epsilon^2}}-1\r).
$$
Using that $y^p={p\over \Gamma(1-p)} \int_0^\infty s^{-p-1}(1-e^{-sy}) ds$, followed by an integration by parts, yields
$$
\E((1+\psi_o/\epsilon)^p)={1\over\Gamma(1-p)}\int_0^\infty
e^{-H(s)}H'(s)s^{-p}\,ds.
$$
The boundary terms vanish since $H(s)=O(s)$ at zero and $s^{-p}\to0$ at infinity. The function $H$ is strictly increasing, and the inverse change of variables $v=H(s)$ is
$$
s=v[1-x\chi(Rv/r)].
$$
As $I_{r,S}/\eta_o^p=\E((1+Y)^p)$, integration against the gamma density of $T$ gives the second formula in \eqref{one_point_common_rep}.
\end{proof}

The role of the constant defined in \eqref{one_point_constants} appears in the identities below (which are also related to the limiting values of $I_A$ and $I_{S,r}$ when $\epsilon\to \infty$, see later).
\begin{eqnarray}\label{one_point_first_moments}
\E(\chi(U))={2a_p\over p},\qquad
\E(\chi(T/r))={2b_p(r)\over p}.
\end{eqnarray}
The first equality follows from integration by part with
$$
{\Gamma(p)\Gamma(1-p)} \E(\chi(U))
=2\int_0^\infty  {u^{-p-1}\over 1+u} (1-\sqrt{1+u})du =\int_0^\infty {u^{-p}\over p}(({1\over 1+u})^{3/2}-2({1\over 1+u})^2) du
$$
and computation of beta-prime integrals.
The second identity follows by inserting
$$
(1+u)^p-1={p\over\Gamma(1-p)}\int_0^\infty
e^{-t}(1-e^{-ut})t^{-p-1}\,dt
$$
in the definition of $b_p$, and using
$$
\int_0^\infty e^{-ru}u^{-3/2}(1-e^{-tu})\,du
=2\sqrt\pi(\sqrt{r+t}-\sqrt r).
$$
Since $\chi$ is strictly decreasing from $1$ to $0$, $\E(\chi(T/r))$  is strictly increasing and continuous from $0$ to $1$ as $r$ goes from $0$ to $\infty$. The first expectation in \eqref{one_point_first_moments} lies strictly between $0$ and $1$. This proves the existence and uniqueness of $r_*(p)$, as well as the assertions following \eqref{one_point_constants}.

At $r=r_*(p)$, the random variables $\chi(U)$ and $\chi(T/r)$ have the same mean. The ratio of the densities of $T/r$ and $U$ is
$$
{f_{T/r}(u)\over f_U(u)}=\Gamma(p)r^{1-p}(1+u)e^{-ru},\qquad
{d\over du}\log{f_{T/r}(u)\over f_U(u)}={1\over1+u}-r.
$$
If $r<1$, this ratio is unimodal and tends to zero at infinity. When $r=r_*(p)$ it must be smaller than $1$ at $u=0$: indeed, otherwise the two densities have a unique crossing point, and since the function $\chi$ is strictly decreasing it prevents the equality $\E(\chi(U))=\E(\chi(T/r_*))$. As both densities integrate to one, their ratio then crosses $1$ exactly twice. 
Since $\chi$ is strictly decreasing it implies that $f_{\chi(U)}-f_{\chi(T/r_*)}$ have signs $+$, $-$, $+$ on successive intervals $[0,y_1]$, $[y_1,y_2]$, $[y_2,1]$, where  $f_{\chi(U)}$ and $f_{\chi(T/r_*)}$ and are the densities of $\chi(U)$ and $\chi(T/r_*)$.  
If $g$ is a convex function and  $\ell$ is the affine function matching with $g$ at points $y_1$ and $y_1$, we have, since $\E\l(\chi(T/r_*)\r)=\E\l(\chi(U)\r)$,
\begin{eqnarray}\label{one_point_convex_order}
\nonumber \E\l(g(\chi(U))\r)-\E\l(g(\chi(T/r_*))\r)&=& \E\l((g-\ell)(\chi(U))\r)-\E\l((g-\ell)(\chi(T/r_*))\r)
\\
&=& \int_0^1 (g(y)-\ell(y)) (f_{\chi(U)}-f_{\chi(T/r_*)}) dy\ge 0,
\end{eqnarray}
since the integrant is point-wise non-negative by convexity. The inequality is strict if $g$ is strictly convex. Hence we proved a convex order $\chi(T/r_*) \leq_{\mathrm{cx}} \chi(U)$.

For $R\ge1$, the transport from the variable $\chi(u)$ to $\chi(Ru)$ is
$$
k_R(y)=\chi(R\chi^{-1}(y))
={2y\over y+\sqrt{y^2+4R(1-y)}},\qquad 0<y<1.
$$
It is increasing and convex. Indeed, writing $k=k_R(y)$, we have
$$
{1-k\over k^2}=R{1-y\over y^2},\qquad 0<k\le y<1,
$$
and implicit differentiation gives
$$
k_R'(y)={k(1-k)(2-y)\over y(1-y)(2-k)}>0,
\qquad
k_R''(y)={2k(1-k)(y-k)(4-k-y)\over
y^2(1-y)^2(2-k)^3}\ge0.
$$
Consequently, $y\mapsto(1-xk_R(y))^{-p}$ is strictly convex for $0<x<1$. Applying \eqref{one_point_convex_order} to this function and then \eqref{one_point_common_rep} proves $I_{r_*,S}<I_A$ for all $a\in \R$ and $\epsilon>0$.

To see that the lower constant is optimal, take $a=0$, so $R=1$, and let $\epsilon\to\infty$, so $x\downarrow0$. The expansion $(1-xy)^{-p}=1+pxy+O_p(x^2)$ is uniform for $0\le y\le1$. By \eqref{one_point_first_moments},
$$
{I_A-I_{r,S}\over\eta_o^p}=2[a_p-b_p(r)]x+O_p(x^2).
$$
For every $r>r_*(p)$ the right-hand side is negative when $\epsilon=1/x-1$ is large enough.

It remains to prove the upper bound. At $r=p$, the density ratio above satisfies
$$
{f_{T/p}(0)\over f_U(0)}=\Gamma(p)p^{1-p}
={\Gamma(1+p)\over p^p}>1.
$$
Indeed, $g(p)=\log\Gamma(1+p)-p\log p$ vanishes at $0$ and $1$, and
$$
g''(p)=\sum_{n=1}^\infty{1\over(n+p)^2}-{1\over p}<0,
\qquad 0<p<1.
$$
Thus $g(p)>0$. The density ratio is unimodal and tends to zero at infinity, so it crosses $1$ exactly once, from above. Consequently $T/p$ is strictly smaller than $U$ in stochastic order. Since $u\mapsto(1-x\chi(Ru))^{-p}$ is strictly decreasing, the common representation \eqref{one_point_common_rep} implies $I_A<I_{p,S}$. This completes the proof.
\end{proof}

\begin{appendix}
\section{Reduction of block-$\mathbb{H}^{2 | 2}$: proof of Lemma~\ref{lem_block}}\label{appendix_A}

Denote by $G_\beta=(\beta-W)^{-1}$ the Green function on $\cV$ and $\tilde G_{\tilde \beta}=(\tilde \beta-\tilde W)^{-1}$ the Green function of the block model. 
Using \eqref{G_paths}, for all $i,j\in \cV$ such that $i\in B_{x}$, $j\in B_{y}$, $x$, $y\in X$,  we have
$$
G_\beta(i,j)=\sum_{\sigma} {\prod_{k=0}^{\vert\sigma\vert-1} W_{\sigma_k,\sigma_{k+1}}\over \prod_{k=0}^{\vert\sigma\vert}\beta_{\sigma_k}}.
$$
where the sum runs on paths $\sigma$ in $\cV$ from $i$ to $j$. Decomposing the path on the successive blocks, we have easily
$$
G_{\beta}(i,j)=\sum_{\tilde \sigma} (\hat G^{B_{x}}  W_{B_{x}, B_{z_1}}  \hat G^{B_{z_1}} W_{B_{z_1}, B_{z_2}} \cdots W_{B_{z_{\vert \tilde \sigma\vert -1}}, B_{y}}  \hat G^{B_{y}})_{i,j}.
$$
where the sum runs on paths $\tilde \sigma=(z_0=x, z_1, \ldots, {z_{\vert \tilde \sigma\vert}}=y)$ on $X$ from $x$ to ${y}$. Using that $W$ is block matrix, hence constant on each term of the form $ W_{B_{z_{k-1}}, B_{z_k}}$, we have
\beq
&&
\nonumber {\left< \indic_{B_{x}}, G \indic_{B_{y}} \right>\over  \vert B_{x}\vert  \vert B_{y}\vert}
\\
\nonumber
&=&
\sum_{\tilde \sigma} {\left< \indic_{B_{x}}, \hat G^{B_{x}} \indic_{B_{x}}\right>\over \vert B_{x}\vert^2}  \tilde W_{x,z_1}  {\left< \indic_{B_{z_1}},\hat G^{B_{z_1}} \indic_{B_{z_1}}\right>\over \vert B_{z_1}\vert^2} 
\tilde W_{z_1, z_2} \cdots \tilde W_{{z_{\vert \tilde \sigma\vert -1}}, {y}}  {\left< \indic_{B_{y}} \hat G^{B_{i_{y}}} \indic_{B_{y}}\right>\over \vert B_{y}\vert^2} 
\\
&=&
\tilde G_{\tilde\beta}(x, y)
\eeq
It in particular implies that
\begin{eqnarray}\label{psi}
\nonumber {1\over \vert B_{x}\vert } \left< \indic_{B_{x}},\psi\right>&=& {1\over \vert B_{x}\vert }  \left< \indic_{B_{x}}, G \eta\right>
\\
\nonumber &=&\sum_{y}  {1\over \vert B_{x}\vert }  \left< \indic_{B_{x}}, G \indic_{B_{y}}  \right>{1\over \vert B_{y}\vert} \tilde \eta_{y}
\\
&=&
(\tilde G_{\tilde\beta}\tilde\eta)_{x}.
\end{eqnarray}

To prove that $(\tilde \beta_x)_{x\in X}$ has law $\nu_X^{\tilde W,\tilde \eta}$ we proceed by recurrence on the size of $X$. Take $x\in X$: we will prove that conditioned on $\beta_{B_{x}^c}$,  $\tilde \beta_{x}$ defined by \ref{beta_tilde} has the law of $\tilde \beta_{{x}}$ conditioned on $\tilde\beta_{X\setminus\{x\}}$ under $\nu_X^{\tilde W,\tilde \eta}$. Then, we will show that the restriction property applied to $\beta_{B_{x}^c}$ implies that $(\tilde \beta_{y})_{y\in X\setminus\{x\}}$ defined by \eqref{beta_tilde} has the law $\nu_X^{\tilde W,\tilde \eta}$ restricted to $X\setminus\{x\}$.

Consider $\beta\sim\nu_\cV^{W,\eta}$, then $\beta_{B_{x}}$ conditioned on $\beta_{B_{x}^c}$ has law $\nu_{B_{x}}^{\check W, \check \eta}$ where the value of $\check W$ is given by Proposition \ref{lem:conditionalbeta} but will not be used here, and
\beq
\check \eta&=&
\eta_{B_{x}}+\sum_{j\in X, j\neq x} W_{B_{x},B_{j}} \hat G^{B_{x}^c}\eta_{B_{x}^c}
\\
&=&
{1\over \vert B_{x}\vert} \indic_{B_{x}} \tilde\eta_{x}+{1\over \vert B_{x}\vert} \indic_{B_{x}} \sum_{j\in X, j\neq x}\tilde W_{x,j} \hat{\tilde G}^{X\setminus\{x\}} \tilde \eta_{X\setminus\{x\}}
\\
&=&
{1\over \vert B_{x}\vert} \indic_{B_{x}} \check{\tilde \eta}_{x},
\eeq
where $\check{\tilde \eta}_{x}$ is the boundary term appearing in the conditioning part of Proposition \ref{lem:conditionalbeta} for the restriction to $\{x\}$ for $\tilde \beta$, $\tilde W$ and $\tilde \eta$ given by \eqref{beta_tilde}, \eqref{W_tilde} and \eqref{eta_tilde}. Hence, by Lemma~\ref{lem_eta-psi} below, $\tilde \beta_{x}$ given by formula \eqref{beta_tilde},
 conditioned on $\beta_{B_{x}^c}$ under $\nu_\cV^{W,\eta}$ has the same law as $\tilde \beta_{x}$ conditioned on $\tilde \beta_{X\setminus\{x\}}$ under $\nu_X^{\tilde W,\tilde \eta}$. Indeed, we have
 \begin{eqnarray}\label{beta_tilde_calcul}
\tilde \beta_{x}=
{\vert B_{x}\vert^2 \over \left<\indic_{B_{x}}, \hat G^{B_{x}} \indic_{B_{x}}\right>}=
{\left<\eta_{B_{x}},\indic_{B_{x}}\right>^2 \over \left<\eta_{B_{x}} , \hat G^{B_{x}} \eta_{B_{x}} \right>}
\end{eqnarray}
since $\eta$ is constant on $B_{x}$ by hypotheses.

To complete the proof we just need to verify that $(\tilde \beta)_{j\in X\setminus\{x\}}$ has the law of $\nu_{X}^{\tilde W,\tilde \eta}$ restricted to $X\setminus\{x\}$. But this a consequence of the restriction property applied to $\nu_\cV^{W,\eta}$ applied to $\cV\setminus B_{x}$ and the recurrence hypothesis.

Note that the initialization is a direct consequence of the definition and Lemma~\ref{lem_eta-psi} below following the computation done above in \eqref{beta_tilde_calcul}.

\begin{lemma}\label{lem_eta-psi}
Let $\cV$ be finite and $W$, $\eta$ be some conductance matrix and boundary term. 
Then,
$$
{\left< \eta, G_\beta \eta\right> \over \left< \eta,1\right>^2}
$$
has inverse Gaussian law with parameters IG$({1\over \left< \eta, 1\right>}, 1)$. In particular it does depend on $W$ and it has the same law as $1\over \beta_1$ for the single point model with boundary $\beta_1\sim \nu_{\{1\}}^{0,  \left< \eta, 1\right>}$.
\end{lemma}
\begin{proof}
Look at the Laplace transform of $\left< \eta, G_\beta \eta\right>$. The term $\exp(-\ t \left< \eta, G_\beta \eta\right>)$ can be combined with the term $\exp(-\demi \left< \eta, G_\beta \eta\right>)$ which appears in the density $\nu_{\cV}^{W, \eta}$ to create the density $\nu_{\cV}^{W,\sqrt{1+2t}\eta}$ but it leaves a constant term $\exp(-(\sqrt{1+2t}-1) \left< \eta, 1\right>)$ equal to the Laplace transform of the law IG$({1\over \left< \eta, 1\right>}, 1)$.
\end{proof}

\end{appendix}
\section*{ Acknowledgement.}
I warmly thank Alain-Sol Sznitman and Xiaolin Zeng for interesting discussions on these questions.

This work is supported by the  project  ANR LOCAL (ANR-22-CE40-0012-02)  operated by the Agence Nationale de la Recherche (ANR), and the Institut Universitaire de France.
\section*{AI declaration}
We started to work on these questions long time ago and a large part of the paper has been done with virtually no AI: this includes Section 1 to 3, where we only use the free version of Gemini in the single Lemma \ref{lem_Green_estimate} to find references on the asymptotic behavior of Green's functions and to understand the proof in \cite{Vainberg2001}. It also includes Section~4 with a more general version of the conjecture.

However, before posting the paper on arXiv, we decided to test Conjecture~\ref{conjecture} with ChatGPT6, both theoretically and numerically. It found a counterexample to our initial  version of the conjecture. It was able to prove the present version in a few limited cases and performed several numerical tests. Section~\ref{sec_one_point} is based on these interactions.

\footnotesize
\bibliographystyle{plain}
\bibliography{AS.bib}  
\end{document}